\documentclass[11pt, a4paper]{article}

\usepackage[ansinew]{inputenc}

\usepackage{amsmath}
\usepackage{amssymb}
\usepackage{amstext}
\usepackage{amsfonts}
\usepackage{enumerate}
\usepackage{tikz}

 \relax

\title{Subvarieties of Pseudocomplemented Kleene Algebras}

\author{\sc D. Casta\~no, V. Casta\~no, J. P. D\'iaz Varela and M. Mu\~noz Santis }

\date{}

\newtheorem{Theorem}{Theorem}[section]
\newtheorem{Proposition}[Theorem]{Proposition}
\newtheorem{Lemma}[Theorem]{Lemma}
\newtheorem{Definition}[Theorem]{Definition}
\newtheorem{Corollary}[Theorem]{Corollary}
\newtheorem{Remark}[Theorem]{Remark}
\newtheorem{Example}[Theorem]{Example}

\newenvironment{Proof}{\noindent \bf Proof. \rm}{$\quad \hfill \blacksquare$}

\usepackage{pictex}

\begin{document}

\maketitle

\begin{abstract}
  In this paper we study the subdirectly irreducible algebras in the
variety  ${\cal PCDM}$ of pseudocomplemented De Morgan algebras by
means of their De Morgan $p$-spaces. We introduce the notion of {\it
body} of an algebra ${\bf L} \in {\cal PCDM}$ and determine $Body({\bf L})$ when
${\bf L}$ is subdirectly irreducible. As a consequence of this, in the case of pseudocomplemented Kleene algebras, three special subvarieties arise naturally, for which we give explicit identities that characterize them. We also introduce a
subvariety ${\cal BPK}$ of ${\cal PCDM}$, namely the variety of {\it bundle
pseudocomplemented Kleene algebras}, determine the whole subvariety lattice and find explicit equational  bases for each of the subvarieties.  In addition, we study the subvariety ${\cal BPK}_0$ of ${\cal BPK}$ generated by the simple members of ${\cal BPK}$, determine the structure of the free algebra over a finite set and their finite weakly projective algebras.

\end{abstract}

\section{Introduction and Preliminaries}

A {\it pseudocomplemented De Morgan algebra} ($pm$-algebra for short) is an algebra 
\linebreak $\langle L, \wedge, \vee, ', *, 0, 1 \rangle$  \ such that $\langle
L, \wedge, \vee, ', 0, 1 \rangle$ is a {\it De Morgan algebra} and
$*$ is a pseudocomplement on $L$, that is, $*$ is a mapping from $ L$
into $L$ such that $a \wedge b = 0 \Leftrightarrow b \leq a^*$. We
denote this variety by ${\cal PCDM}$.

\medskip

In 1981, Romanowska \cite{Roma} initiated the study of the variety
of De Morgan algebras with  pseudocomplementation by characterizing
its finite subdirectly irreducible members. In 1986, Sankappanavar
\cite{San2} began an investigation of a larger variety, namely the
variety of Ockham  algebras with pseudocomplementation. In
particular, he characterized the subdirectly irreducible, non
regular, pseudocomplemented De Morgan algebras extending
Romanowska's results. Sankappanavar's results motivated us to
attempt to characterize all subdirectly irreducible
pseudocomplemented De Morgan algebras.

Recently in \cite{San2020} H. Sankapannavar,  M. Adams and J. Vaz de Carvalho studied some subvarieties of pseucomplemented De Morgan algebras; more pecisely the subvarieties ${\bf M}_n$, resp. ${\bf K}_n$, of regular De Morgan, resp. Kleene, pseudocomplemented algebras of rank $n$. In that article they fully described the lattice of subvarieties for particular values of $n$. With the intention of furthur studying the lattice of subvarieties of pseudocomplemented De Morgan and Kleene algebras we introduced the variety of {\it bundle pseudocomplemented Kleene algebras} (not contained in any of the varieties studied in \cite{San2020}) and fully described its subvariety lattice, giving all equational bases for each subvariety. We also studied the subvariety ${\cal BPK}_0$, which coincides with the variety ${\bf K}_1$ in \cite{San2020}, in depth. In addition to giving equational bases for each subvariety, we characterized free algebras on a finite number of generators and finite weakly projective algebras in this variety.

\medskip

The main tool we use in this paper is a duality between the category
of pseudocomplemented De Morgan algebras and certain topological
spaces based on the duality developed by Priestley. Next we give a
brief summary of the necessary required facts; for further
information (see \cite{Pri1}, \cite{Pri2} and \cite{Pri3}).

\medskip

Given a subset $Y$ of a poset $\langle  X, \leq \rangle$, define $[Y)= \{a \in X :
a \geq b \mbox{ for some } b \in Y\}$ and $(Y] = \{a \in X : a \leq
b \mbox{ for some } b \in Y\}$. $Y$ is {\it decreasing} if $P=(P]$ and $P$ is {\it increasing} if $P=[P)$.
\medskip

A triple $\langle  X, \leq, \tau \rangle$ is a {\it totally order disconnected
topological space} if $\langle X, \leq \rangle$ is a poset, $\tau$ is a topology
on $X$, and for $a$, $b \in X$, if $a \not\leq b$ then there exists
a clopen (open and closed)  increasing $V \subseteq X$ such that $a \in V$ and $b
\notin V$. A {\it compact} totally order disconnected space is
called a {\it Priestley space}.
\medskip

In \cite{Pri1} and \cite{Pri3}, H.~A.~Priestley showed that the
category of bounded distributive lattices and $(0,1)$--lattice
homomorphisms is dually equivalent to the category of Priestley
spaces and order preserving continuous functions (see also the
survey paper \cite{Pri2}).
\medskip

If ${\bf X}$ is a Priestley space, then $\langle \mathbb{D}({\bf X}), \cap,
\cup, \emptyset, X \rangle$ 
is the lattice of clopen increasing
subsets of ${\bf X}$. If $f:{\bf X}_1 \to {\bf X}_2$ is a continuous order preserving map,
then $\mathbb{D}(f):\mathbb{D}({\bf X}_2)\to \mathbb{D}({\bf X}_1)$ defined by
$\mathbb{D}(f)(V) = f^{-1}(V)$ is a $(0,1)$--lattice homomorphism.
Conversely, if ${\bf L}$ is a bounded distributive lattice, then the set
of prime filters of ${\bf L}$ denoted $\mathbb{X}({\bf L})$ is a Priestley
space, ordered by set inclusion and with the topology having as a
sub-basis the sets $\eta_{{\bf L}}(a)=\{P \in \mathbb{X}({\bf L}) : a \in P\}$
and $\mathbb{X}({\bf L}) \setminus \eta_{{\bf L}}(a)$ for $a \in L$. If $h:{\bf L}_1\to
{\bf L}_2$ is a $(0,1)$--lattice homomorphism, then
$\mathbb{X}(h):\mathbb{X}({\bf L}_2)\to \mathbb{X}({\bf L}_1)$ defined by
$\mathbb{X}(h)(P)=h^{-1}(P)$ is a continuous order preserving map.
In addition, the mapping $\eta_{{\bf L}}:{\bf L}\to \mathbb{D}(\mathbb{X}({\bf L}))$
is a lattice isomorphism, and $\varepsilon_{{\bf X}}:{\bf X}\to
\mathbb{X}(\mathbb{D}({\bf X}))$ defined by $\varepsilon_{{\bf X}}(x)=\{V \in
\mathbb{D}({\bf X}): x \in V\}$ is a homeomorphism and an order
isomorphism.
\medskip

Since $p$-algebras (distributive lattices with pseudocomplementation) are bounded distributive lattices, the category
of $p$-algebras is isomorphic to a subcategory of bounded
distributive lattices. A {\it p-space} is a Priestley space $\langle X,
\leq, \tau \rangle$ such that $(Y]$ is open for every  $Y \in
\mathbb{D}({\bf X})$. If ${\bf X}_1$ and ${\bf X}_2$ are $p$-spaces, a $p$-morphism is a
continuous order-preserving map $\phi:{\bf X}_1\to {\bf X}_2$ for which $\phi([P) \cap Max({\bf X}_1)) = [\phi(P)) \cap Max({\bf X}_2)$, where $Max({\bf X}_1)$ and $Max({\bf X}_2)$ denote the set of all maximal points of ${\bf X}_1$ and ${\bf X}_2$ respectively.
For a $p$-algebra ${\bf L}$ and $a \in L$, $\eta(a)
= V_a $ denotes the clopen increasing set that represents $a$. If $a
\in L$ then, under the duality given above, $a^\ast$ corresponds to
the clopen increasing set $(V_a]^c$, where $Y^c$ denotes the
complement of $Y$ relative to $\mathbb{ X}({\bf L})$. As a consequence, Priestley's duality leads us to
the following fact: the functors $\mathbb{ X}$ and $\mathbb{D}$
establish a dual equivalence between the category of $p$-algebras
and the category of $p$-spaces (see \cite{Pri2}).
\medskip

If ${\bf X}$ is a Priestley space and $\phi: {\bf X} \to {\bf X}$ is an order-reversing involutive ($\phi = \phi^{-1}$) homeomorphism then
$\langle {\bf X}, \phi \rangle$ is called a {\it De Morgan space} \cite{Bly-Var}. If
$\langle M, \wedge, \vee, ^{\prime}, 0,1 \rangle$ is a De Morgan
algebra and $\phi : \mathbb{X}({\bf M}) \to \mathbb{X}({\bf M})$ is given by
$\phi(P) = P^{\prime c}$, where $P'= \{a' \in M: a \in P\}$, then
$\langle \mathbb{X}({\bf M}), \phi \rangle$ is a De Morgan space. $\phi$ is called
the {\it Birula-Rasiowa transformation}. If $\eta(a)=V_a$ denotes the
clopen increasing set that represents $a \in M$, then under the
duality given above, $a^{\prime}$ corresponds to the clopen
increasing set $\phi (V_a)^c = \mathbb{X}({\bf M}) \setminus \phi
(V_a)$. Conversely, if $\langle {\bf X}, \phi \rangle$ is a De Morgan space, $\langle
\mathbb{D}({\bf X}), \cap, \cup, ^{\prime}, \emptyset,X \rangle$ is a De
Morgan algebra where for $V \in \mathbb{D}({\bf X})$, $V^{\prime} =
\phi(V)^c$. The category whose objects are De Morgan spaces and whose morphisms are
continuous order-preserving functions $f:{\bf X}_1\to {\bf X}_2$ such that $f \circ \phi_1 =
\phi_2 \circ f$ is dually equivalent to the category
of De Morgan algebras and (De Morgan) homomorphisms.
\medskip

A {\it De Morgan $p$-space} is a system $\langle X, \leq, \tau, \phi\rangle$
which is both a $p$-space and a De Morgan space. Morphisms in the
category of De Morgan $p$-spaces will be functions $f: {\bf X} \to
{\bf X}^{\prime}$ which are morphisms in the category of $p$-spaces and in
the category of De Morgan spaces.

It is immediate that the category of De Morgan $p$-spaces and the
category of pseudocomplemented De Morgan algebras are dually equivalent.

\medskip

\section{Simple and subdirectly irreducible pseudocomplemented De Morgan algebras}

It is known that for each $a$ in any Priestley space ${\bf X}$, there
exists at least one point $b \geq a$ such that $b$ is maximal with
respect to the partial order. The set of all maximal points of ${\bf X}$
will be denoted by $Max({\bf X})$. We will write it simply $Max$ when no confusion may arise.
\medskip

A subset $Y$ of a De Morgan space $\langle {\bf X}, \phi \rangle$ is called an {\it
involutive subset} if it is $\phi$-invariant, in the sense that
$\phi(a) \in Y$ whenever $a \in Y$.
\medskip

It is well known (see, for instance, \cite[Th. 4.2]{Bly-Var}) that
there exists a one-to-one correspondence between the congruences on
a De Morgan algebra ${\bf L}$ and the closed and involutive subsets of
$\mathbb{ X}({\bf L})$. The correspondence is given by $Y \longmapsto
\theta(Y)$,  where $$(U,V) \in \theta(Y) \Longleftrightarrow U \cap
Y = V \cap Y.$$
\medskip

On the other hand there exists a one-to-one correspondence between
the congruences on a $p$-algebra ${\bf L}$ and
the closed subsets $Y$ of $\mathbb{X}({\bf L})$ satisfying
$Max(\mathbb{ X}({\bf L})) \cap [Y) \subseteq Y$ (\cite[Proposition
2.3]{Dav}).
\medskip

From this we have the following result.

\begin{Proposition} \label{congruencias}
Let ${\bf L}$ be a $pm$-algebra. The lattice of
congruences on ${\bf L}$ is dually isomorphic to the lattice of closed and
involutive sets $Y$ in $\mathbb{X}({\bf L})$ that satisfy
$Max(\mathbb{X}({\bf L})) \cap [Y) \subseteq Y$. This anti-isomorphism
$\theta$ is given by $Y \longmapsto \theta(Y)$,  where $$(U,V) \in
\theta(Y) \Longleftrightarrow U \cap Y = V \cap Y.$$

\end{Proposition}

Let ${\bf L} \in {\cal PCDM}$. A subset $Y$ of $\mathbb{X}({\bf L})$ is called a {\it C-subset} \ if $Y$
is closed, involutive and satisfies $Max(\mathbb{X}({\bf L})) \cap [Y)
\subseteq Y$. It is easy to see that if $Y$ is a $C$-subset, then $Y$
also satisfies $Min(\mathbb{X}({\bf L})) \cap (Y] \subseteq Y$, where
$Min(\mathbb{X}({\bf L}))$ denotes the set of all minimal points of
$\mathbb{X}({\bf L})$ ($Min$ for short).
\medskip

The following result shows that $Max(\mathbb{X}({\bf L})) \cup Min(\mathbb{ X}({\bf L}))$ is a $C$-subset  for all $pm$-algebras.

\medskip

\begin{Lemma} \label{max y min es cerrado}
$Max(\mathbb{X}({\bf L})) \cup Min(\mathbb{X}({\bf L}))$ is  closed in
$\mathbb{X}({\bf L})$, for every pseudocomplemented De Morgan algebra ${\bf L}$.
\end{Lemma}
\begin{Proof}
It is known that $Max$ is a closed set in
$\mathbb{X}({\bf L})$. Indeed, if $P \in \mathbb{X}({\bf L}) \setminus
Max$, take $U \in Max$ such that $P
\subseteq U$. Then $U \not \subseteq P$, and consequently there exists a clopen
increasing set $V$ with $U \in V$ and $P\notin V$. Then $(V]$ is
open and so $(V] \setminus V$ is an open naighborhood of $P$
disjoint of $Max$.

Since $\phi(Max)=Min$ and $\phi$ is a
homeomorphism, then $Min$ is a closed set, and so is
$Max \cup Min$.
\end{Proof}

\medskip

A {\it connected component} of a partially ordered set ${\bf X}$ is a
non-empty subset of ${\bf X}$ that is both increasing and decreasing and
it is minimal with respect to these properties. ${\bf X}$ is the disjoint
union of its connected components. ${\bf X}$ is said to be {\it connected}
if ${\bf X}$ has exactly one connected component.
\medskip

Let $\mathbb{X}({\bf L}) = \bigcup_{i \in I} S_i$, where $\{S_i\}_{i \in
I}$ is the collection of connected components of $\mathbb{X}({\bf L})$. If
$\phi(S_i) \subseteq S_j$, then it is easily seen that $\phi(S_i) =
S_j$ and $\phi(S_j) = S_i$, and consequently, $\phi(S_i \cup S_j) =
S_i \cup S_j$. In that case we say that $S_i \cup S_j$ is a
$\phi$-{\it connected component} of $\mathbb{X}({\bf L})$. A $\phi$-{\it
connected component} of $\mathbb{X}({\bf L})$ can be defined as a
non-empty subset of $\mathbb{X}({\bf L})$ that is increasing, decreasing
and involutive and it is minimal with respect to these properties.

We say that $\mathbb{X}({\bf L})$ is $\phi$-{\it connected} if
$\mathbb{X}({\bf L})$ has exactly one $\phi$-connected component. In this case we have that $Max \cup Min \subseteq Y$ for each nonempty $C$-subset $Y$ as it shows in the following lemma.

\medskip

\begin{Lemma} \label{C-subconjuntos en phi-conexos}
Let ${\bf L}$ be a $pm$-algebra such that
$\mathbb{X}({\bf L})$ is $\phi$-connected. If $Y$ is a $C$-subset, $Y \neq
\emptyset$, then $Max(\mathbb{X}({\bf L})) \cup Min(\mathbb{X}({\bf L}))
\subseteq Y$.
\end{Lemma}
\begin{Proof}
In order to prove that $Max \cup Min
\subseteq Y$ we will show that $(Max \cap Y] =
[Min \cap Y)$.

If $P \in (Max \cap Y]$, $P \in (Y]$. Let $M$ be a minimal filter such that $M \subseteq P$, so $M \in
Min \cap (Y]$. Since $Y$ is a $C$-subset, we have that
$M \in
Min \cap (Y] \subseteq Y$ and so $M \in Y$. Then $P \in [Y \cap Min)$. In the same way we can see that $[Min \cap Y) \subseteq (Max \cap Y]$.

From this we have that $(Max \cap Y] =
[Min \cap Y)$ is increasing, decreasing and
involutive. Since $\mathbb{X}({\bf L})$ is $\phi$-connected we have that
$(Max \cap Y] = [Min \cap Y)=
\mathbb{X}({\bf L})$. Therefore $Max \cup
Min \subseteq Y$.
\end{Proof}

\medskip

The following results characterize the simple and directly
indecomposable algebras  ${\bf L} \in {\cal PCDM}$ by means of their dual
space, where $\mathbb{X}({\bf L})$ is $\phi$-connected.

\medskip

\begin{Theorem} \label{simple}

Let ${\bf L}$ be a $pm$-algebra such that
$\mathbb{X}({\bf L})$ is $\phi$-connected. Then ${\bf L}$ is simple if and only
if $\mathbb{X}({\bf L})=Max(\mathbb{X}({\bf L})) \cup Min(\mathbb{X}({\bf L}))$.
\end{Theorem}

\begin{Proof}
Suppose that $\mathbb{X}({\bf L}) \ne Max \cup
Min$. Then the congruence $\theta_Y$ associated to
the $C$-subset $Y=Max \cup Min$ is not
a trivial congruence, and so ${\bf L}$ is not simple.

Conversely,   suppose that $\mathbb{X}({\bf L})$ is connected and $\mathbb{X}({\bf L})= Max \cup Min$. If  $Y \ne \emptyset$ is a $C$-subset, by Lemma \ref{C-subconjuntos en phi-conexos}, we have that $Max \cup
Min \subseteq Y$. Consequently,  $Y=\mathbb{ X}({\bf L})$, and so ${\bf L}$ is simple.

\end{Proof}

\medskip

\begin{Example} {\rm \label{si infinita pero no simple}

Let us consider the set $ X_1= \mathbb{N} \cup \{\infty\}$ with the topology given by the  one-point compactification of a countable discrete space topology and the set $X_2=\{0,1\}$ with the discrete topology. The following picture gives an order $\leq$ for $ X_1 \times  X_2=A$. For simplicity of notation, we use $n$  for $(n, 1)$ and $n'$ for $ (n,0)$, with $n\in\mathbb{N}\cup \{\infty\}$.

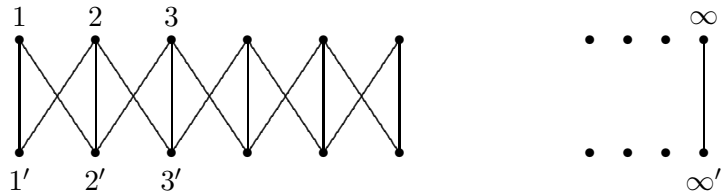
$$
\beginpicture
\setcoordinatesystem units <0.50000cm,0.50000cm>
\setlinear
\setshadesymbol <z,z,z,z> ({\fiverm .})
\putrule from -9.00000 0.00000 to -9.00000 -3.00000
\putrule from -7.00000 0.00000 to -7.00000 -3.00000
\putrule from -5.00000 0.00000 to -5.00000 -3.00000
\putrule from -3.00000 0.00000 to -3.00000 -3.00000
\putrule from -1.00000 0.00000 to -1.00000 -3.00000
\putrule from 1.00000 0.00000 to 1.00000 -3.00000
\putrule from 9.00000 0.00000 to 9.00000 -3.00000
\plot -9.00000 -3.00000 -7.00000 0.00000 /
\plot -9.00000 0.00000 -7.00000 -3.00000 /
\plot -7.00000 -3.00000 -5.00000 0.00000 /
\plot -7.00000 0.00000 -5.00000 -3.00000 /
\plot -5.00000 0.00000 -3.00000 -3.00000 /
\plot -5.00000 -3.00000 -3.00000 0.00000 /
\plot -3.00000 0.00000 -1.00000 -3.00000 /
\plot -3.00000 -3.00000 -1.00000 0.00000 /
\plot -1.00000 0.00000 1.00000 -3.00000 /
\plot -1.00000 -3.00000 1.00000 0.00000 /
\put {$\scriptstyle\bullet$} at -9.00000 -3.00000
\put {$\scriptstyle\bullet$} at -9.00000 0.00000
\put {$\scriptstyle\bullet$} at -7.00000 -3.00000
\put {$\scriptstyle\bullet$} at -7.00000 0.00000
\put {$\scriptstyle\bullet$} at -5.00000 -3.00000
\put {$\scriptstyle\bullet$} at -3.00000 -3.00000
\put {$\scriptstyle\bullet$} at -1.00000 -3.00000
\put {$\scriptstyle\bullet$} at 1.00000 -3.00000
\put {$\scriptstyle\bullet$} at 9.00000 -3.00000
\put {$\scriptstyle\bullet$} at -5.00000 0.00000
\put {$\scriptstyle\bullet$} at -3.00000 0.00000
\put {$\scriptstyle\bullet$} at -1.00000 0.00000
\put {$\scriptstyle\bullet$} at 1.00000 0.00000
\put {$\scriptstyle\bullet$} at 9.00000 0.00000
\put {$\scriptstyle\bullet$} at 6.00000 -3.00000
\put {$\scriptstyle\bullet$} at 7.00000 -3.00000
\put {$\scriptstyle\bullet$} at 8.00000 -3.00000
\put {$\scriptstyle\bullet$} at 6.00000 0.00000
\put {$\scriptstyle\bullet$} at 7.00000 0.00000
\put {$\scriptstyle\bullet$} at 8.00000 0.00000
\put {$ 1'$} at -9 -3.7
\put {$ 2'$} at -7 -3.7
\put {$3'$} at -5 -3.7
\put {$ 1$} at -9 0.6
\put {$2$} at -7 0.6
\put {$3$} at -5 0.6
\put {$\infty$} at 9 0.6
\put {$\infty'$} at 9 -3.7
\endpicture
$$

\medskip

\bigskip

If we define the order-reversing involutive homeomorphism $g:A\to A$ such that $g(n)=n'$ and  $g(n')=n$, it is easy to check that $\langle A, \leq, \tau, g \rangle$ is a De Morgan $p$-space. Taking into account that $\{\infty, \infty'\}$  is a closed set but not an open set, we have that $\mathbb{D}(A)$ has a unique non trivial congruence $\theta$, associated  with the $C$-subset  $\{\infty, \infty'\}$.

\medskip

}
\end{Example}

The previous example shows that there exist De Morgan $p$-algebras ${\bf L}$ such that $\mathbb{X}({\bf L}) = Max \cup Min$ and ${\bf L}$ is not a simple algebra.
Moreover, it also shows that the ${\cal PCDM}$ is not locally finite. It suffices to consider the clopen increasing set $V=\{1\}$ and apply $\ast\prime$ successively.

\medskip

In order to study the subdirectly irreducible
pseudocomplemented De Morgan algebras we call {\it body of
${\bf L}$} to $Body({\bf L})=\mathbb{X}({\bf L})\setminus
(Max(\mathbb{X}({\bf L}))\cup Min (\mathbb{X}({\bf L})))$. The following result
gives a necessary condition for an algebra ${\bf L} \in {\cal PCDM} $ to
be  subdirectly irreducible.

We will denote by $|X|$ the number of elements of $X$.

\medskip

\begin{Proposition} \label{si cardinal mayor a 3}
If ${\bf L}$ is a non-trivial subdirectly irreducible $pm$-algebra then  \ $|Body({\bf L})| \leq 2$.
\end{Proposition}
\begin{Proof}
Let us consider $Body({\bf L})= \bigcup_{i \in I} \{P_i,
\phi(P_i)\}$ where $|I| \geq 2$ and the subsets $Y_i=Max \cup
Min \cup \{P_i,\phi(P_i)\}$ for each $i \in I$. By
Lemma  \ref{max y min es cerrado} we have
 that $Y_i$ is a $C$-subset
for each $i \in I$. Then, since $\bigcup_{i \in
I}Y_i=\mathbb{X}({\bf L})$, $\bigcap_{i \in I} \theta_{Y_i}=\Delta$
where $\theta_{Y_i}\ne \Delta$ for all $i \in I$. From this ${\bf L}$
is not a subdirectly irreducible algebra.
\end{Proof}

\medskip

The proof of the proposition may be used to show a stronger result, namely, in the body of a subdirectly irreducible $pm$-algebra there are not  incomparable elements such that $P= \phi(P)$ and  $Q=\phi(Q)$.

\begin{Corollary} \label{body de s.i.}
Let ${\bf L}$ be a $pm$-algebra. If ${\bf L}$ is
subdirectly irreducible, then $|Body({\bf L})| <2$ or
$Body({\bf L}) = \{P, \phi(P)\}$, with $\phi(P) \ne P$.
\end{Corollary}

\begin{Theorem} \label{body y phi-conexo}
Let ${\bf L}$ be a $pm$-algebra such that
$\mathbb{X}({\bf L})$ is $\phi$-connected. If $|Body({\bf L})| <2$ or
$Body({\bf L}) = \{P, \phi(P)\}$, with $\phi(P) \ne P$ then
${\bf L}$ is subdirectly irreducible.
\end{Theorem}

\begin{Proof}
If $|Body({\bf L})| =0$, by Theorem \ref{simple}, ${\bf L}$ is
simple, so ${\bf L}$ is subdirectly irreducible. If
$|Body({\bf L})|=1$ and $Y$ is a $C$-subset, since
$\mathbb{X}({\bf L})$ is $\phi$-connected, by Lemma \ref{C-subconjuntos en
phi-conexos} we have that $Max \cup
Min \subseteq Y$. So $Y=Max \cup
Min$ or $Y=\mathbb{X}({\bf L})$, and then ${\bf L}$ is
subdirectly irreducible. Finally, if $Body({\bf L}) = \{P,
\phi(P)\}$, with $\phi(P) \ne P$, then it is easy to see that the
only non-empty $C$-subsets of $\mathbb{X}({\bf L})$ are $Y_1=
Max \cup Min$ and $Y_2=
\mathbb{X}({\bf L})$. So ${\bf L}$ is subdirectly irreducible.
\end{Proof}
\medskip

Note that the Example \ref{si infinita pero no simple} shows that  there exists subdirectly irreducible $pm$-algebras whose dual spaces associated are not $\phi$-connected.

\medskip

Since  for each  finite subdirectly irreducible algebra ${\bf L}$ in
${\cal PCDM}$, $\mathbb{X}({\bf L})$ is $\phi$-connected, the above
theorems provide necessary and sufficient conditions for a finite pseudocomplemented De
Morgan algebra to be subdirectly irreducible.

\begin{Corollary} \label{caracterizacion s.i. finita}
Let ${\bf L}$ be a finite $pm$-algebra. Then ${\bf L}$ is subdirectly
irreducible if and only if $|Body({\bf L})| <2$ or
$Body({\bf L}) = \{P, \phi(P)\}$, with $\phi(P) \ne P$.
\end{Corollary}


The following table summarizes the results proved so far for a $pm$-algebra ${\bf L}$:

\

\noindent \begin{tabular}{|c|c|}
\hline
If $\mathbb{X}({\bf L})$ is $\phi$-connected: & {\bf L} is simple $\Longleftrightarrow $ $ \mathbb{X}({\bf L})=Max \cup Min$ \\
\hline
If $\mathbb{X}({\bf L})$ is $\phi$-connected: & {\bf L} is s.i  $\Longleftrightarrow $ $ |Body({\bf L})| < 2 \ \ {\mbox  or  } \ \ Body({\bf L}) = \{P, \phi(P)\}$, $P \ne \phi(P)$ \\
\hline
If ${\bf L}$ is finite: & ${\bf L}$ is s.i.  $\Longleftrightarrow $ $ |Body({\bf L})| < 2 \ \ {\mbox  or  } \ \ Body({\bf L}) = \{P, \phi(P)\}$, $P \ne \phi(P)$ \\
\hline
\end{tabular}

\section{Subdirectly irreducible $pk$-algebras}

An important subvariety of the $pm$-algebras is the pseudocomplemented Kleene algebras ($pk$-algebra for short). These algebras are $pm$-algebras satisfying $x \wedge x^\prime \leq y \vee y^\prime$. In this section we determine the dual spaces associated with a subdirectly irreducible $pk$-algebra along with some of their properties. This will lead us to consider two subvarieties of  $pk$-algebras and  we will give equational bases for each of them.

\medskip

Recall that a De Morgan algebra is a Kleene algebra if and only if
for every $P \in \mathbb{X}({\bf L})$, either $P \subseteq \phi(P)$ or
$\phi(P) \subseteq P$ (see \cite{Mon} for more details). Therefore from Theorem \ref{body de s.i.}, the
following result is immediate.

\begin{Remark} {\rm
If ${\bf L}$ is a subdirectly irreducible pseudocomplemented Kleene
algebra, then
\begin{enumerate}

\item $\mathbb{X}({\bf L}) = Max \cup
Min$, or

\item $\mathbb{X}({\bf L}) = Max  \cup
Min \cup \{P\}$, or

\item $\mathbb{X}({\bf L}) = Max  \cup
Min \cup \{P, \phi(P)\}$, with $\phi(P)
\varsubsetneq P$.
\end{enumerate}}
\end{Remark}

The following lemmas provide an important tool for the rest of the
paper. $P$ denotes an element of
$Body({\bf L})$ such that $\phi(P) \subseteq P$. Recall that
if $\langle {\bf X}, \phi \rangle$ is a De Morgan space, $\langle \mathbb{D}({\bf X}), \cap,
\cup, ', \emptyset, X \rangle$ is a De Morgan algebra, where for $V
\in \mathbb{D}({\bf X})$, $V' = \phi(V)^c$.

\medskip

\begin{Lemma} \label{condiciones de V y V prima}
Let ${\bf L}$ be a  subdirectly irreducible $pk$-algebra and  $V \in \mathbb{D}(\mathbb{X}({\bf L}))$.

\begin{enumerate}
\item[$(1)$] \label{V inters V'} $Q \in V \cap V^\prime$ if only if $Q \in V$ and $\phi(Q) \not \in V$.

\item[$(2)$] \label{V interseccion V' no tiene minimales} $V\cap
V'\subseteq Max \cup \{P\}$.
\end{enumerate}

\end{Lemma}

\begin{Proof}
$(1)$  is immediate from the definition of the operation $\prime$ and the properties of $\phi$.

In order to prove $(2)$ let us consider just the case in which $\mathbb{X}({\bf L}) =
Max  \cup Min \cup \{P,
\phi(P)\}$, with $\phi(P) \varsubsetneq P$. If $Q \in
Min$ then $Q = \phi(U)$, for $U \in
Max$. If we suppose that $\phi(U)\in V\cap V'$, then
$\phi(U)\in V$. Therefore, by $(1)$, $U \not \in V$. However, since ${\bf L}$ is a Kleene algebra, $\phi(U) \subseteq U$, which contradicts the fact that $V$ is an increasing set.

If we suppose that $Q = \phi(P)\in V\cap V'$ then $\phi(P)\in V$ and
$P \not \in V$, which is impossible. \end{Proof}

Observe that $(1)$ is even true for those algebras which are not necessarily subdirectly irreducible.

Taking into account the number of the elements for the $Body({\bf L})$ in a subdirectly irreducible $pk$-algebra we can consider two proper subvarieties of the $pk$-algebras namely: ${\cal PK}_0$, the subvariety generated by the subdirectly irreducible algebras ${\bf L}$ such that $|Body({\bf L})|=0$ and ${\cal PK}_1$, the subvariety generated by the subdirectly irreducible algebras ${\bf L}$ such that $|Body({\bf L})|\leq 1$.

\medskip

Next, we will consider the term $C(x)=(x \wedge x') \vee (x \wedge x')^\ast$ and we will use it to find equtional bases for the subvarieties defined above.

Note that for every $pk$-algebra ${\bf L}$, $C(a)$ is a dense element ($C(a)^\ast=0$). In terms of duality, this means that $Max \subseteq C(V)$, for  each clopen increasing  $V$ that correspond with $a$.

\begin{Theorem} \label{identidades para PK0 y PK1}
Let $C(x)$ be the term defined above. If ${\bf L}$ is a subdirectly irreducible $pk$-algebra:

\begin{enumerate}
\item[$(1)$] ${\bf L}$ satisfies $C(x)^\prime \leq C(x)$, for each $x \in {\bf L}$, if only if $|Body({\bf L})|=0$.
\item[$(2)$] ${\bf L}$ satisfies   $C(x)^\prime \wedge C(x) \leq C(y)$, for each $x \in {\bf L}$ if only if $|Body({\bf L})|\leq 1$.

\end{enumerate}
\end{Theorem}

\begin{Proof}
In order to prove $(1)$, suppose that there exists $P \in Body({\bf L})$. Without loss of generality we can assume $\phi(P) \subseteq P$. We know that there exists a maximal element $U$ such that $P \subseteq U$, so, since $\mathbb{X}({\bf L})$ is a totally order disconnected topological space, there exists a clopen increasing $V_U$ such that $U \in V_U$ and $P \not \in V_U$. We claim that $P \in C(V_U)^\prime$ and $P \not \in C(V_U)$. Indeed, from $P \not \in V_U$, $\phi(U) \not \in V_U$, so by Lemma \ref{V inters V'}, $U \in V_U \cap V_U^\prime$ and consequently $P \not \in (V_U \cap V_U^\prime)^\ast$. Further, from $P \not \in V_U \cap V_U^\prime$ we have that $P \not \in C(V_U)$. Moreover, from $P, \phi(P) \not \in C(V_U)$, $P \in C(V_U)^\prime$. This shows that ${\bf L}$ does not satisfy the equation $C(x)^\prime \leq C(x)$.
Finally, taking into account that $Max \subseteq C(V)$ (Lemma \ref{V interseccion V' no tiene minimales}), if $|Body({\bf L})|=0$, it is immediate that $C(V)^\prime \subseteq C(V)$, for every $V \in \mathbb{D}(\mathbb{X}({\bf L}))$.

In order to prove $(2)$, suppose that $Body({\bf L})=\{\phi(P), P\}$ where $\phi(P) \varsubsetneq P$. Since $\mathbb{X}({\bf L})$ is a totally order disconnected topological space, there exists a clopen increasing $V_P$ such that $P \in V_P$ and $\phi(P) \not \in V_P$. We can show, in the same manner as (1), that $P \in C(V_P) \cap C(V_P)'$.

On the other hand, there exists a maximal element $U$ such that $P \subseteq U$, hence, there exists a clopen increasing $V_U$ such that $U \in V_U$ and $P \not \in V_U$. It is easy to check that $P \not \in C(V_U)$. We have thus shown  that $C(x) \wedge C(x)' \leq C(y)$ does not hold in ${\bf L}$.
\end{Proof}

\section{Bundle Pseudocomplemented Kleene Algebras (${\cal BPK}$)}

In this section we focus our attention in some particular De Morgan $pm$-spaces and their corresponding algebras. Let us consider the subdirectly irreducible pseudocomplemented Kleene algebras
whose dual spaces satisfy the following condition: every  maximal element contains all  non-maximal elements. This is shown in the following pictures.

\bigskip

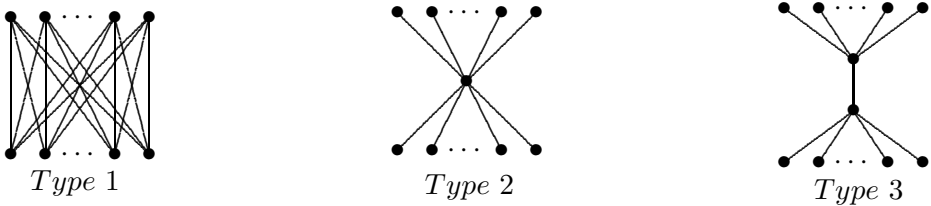
\begin{minipage}{5cm}
\hspace{3cm} \beginpicture \setcoordinatesystem units
<0.4556cm,0.45556cm> \setlinear \setshadesymbol <z,z,z,z> ({\fiverm
.})

\putrule from -4.00000 1 to -4.00000 -3.00000

\plot -4.00000 -3.00000 -3.00000 1 /

\putrule from -3.00000 -3.00000 to -3.00000 1

\plot -3.00000 -3.00000 -1.00000 1 /

\putrule from -1.00000 1 to -1.00000 -3.00000

\plot -1.00000 -3.00000 0.00000 1 /

\putrule from 0.00000 1 to 0.00000 -3.00000

\plot 0.00000 -3.00000 -1.00000 1 /

\plot 0.00000 -3.00000 -3.00000 1 /

\plot 0.00000 -3.00000 -4.00000 1 /

\plot -1.00000 -3.00000 -3.00000 1 /

\plot -1.00000 -3.00000 -4.00000 1 /

\plot -3.00000 -3.00000 -4.00000 1 /

\plot -3.00000 -3.00000 0.00000 1 /

\plot -4.00000 -3.00000 -1.00000 1 /

\plot -4.00000 -3.00000 0.00000 1 /

\put {$\bullet$} at -4.00000 -3.00000

\put {$\bullet$} at -3.00000 -3.00000

\put {$\bullet$} at -1.00000 -3.00000

\put {$\bullet$} at 0.00000 -3.00000

\put {$\bullet$} at -4.00000 1

\put {$\bullet$} at -3.00000 1

\put {$\bullet$} at -1.00000 1

\put {$\bullet$} at 0.00000 1

\put{$\cdots$} at -2 -3

\put{$\cdots$} at -2 1

\put {$Type \ 1$} at -2.13626 -3.82041
\endpicture
\end{minipage}
\begin{minipage}{5cm}
\hspace{3cm}
 \beginpicture \setcoordinatesystem units
<0.45556cm,0.45556cm> \setlinear \setshadesymbol <z,z,z,z> ({\fiverm
.}) \plot -4.00000 -1.00000 -6.00000 1.00000 / \plot -4.00000
-1.00000 -5.00000 1.00000 / \plot -4.00000 -1.00000 -3.00000 1.00000
/ \plot -4.00000 -1.00000 -2.00000 1.00000 / \plot -2.00000 -3.00000
-4.00000 -1.00000 / \plot -3.00000 -3.00000 -4.00000 -1.00000 /
\plot -4.00000 -1.00000 -5.00000 -3.00000 / \plot -6.00000 -3.00000
-4.00000 -1.00000 /

\put {$\bullet$} at -6.00000 -3.00000

\put {$\bullet$} at -5.00000 -3.00000

\put {$\bullet$} at -3.00000 -3.00000

\put {$\bullet$} at -2.00000 -3.00000

\put {$\bullet$} at -4.00000 -1.00000

\put {$\bullet$} at -6.00000 1.00000

\put {$\bullet$} at -5.00000 1.00000

\put {$\bullet$} at -3.00000 1.00000

\put {$\bullet$} at -2.00000 1.00000

\put{$\cdots$} at -4 -3

\put{$\cdots$} at -4 1

\put {$Type \ 2$} at -3.92214 -4.06907
\endpicture
\end{minipage}
\begin{minipage}{5cm}
\hspace{3cm}
\beginpicture
\setcoordinatesystem units <0.45556cm,0.45556cm>  \setlinear
\setshadesymbol <z,z,z,z> ({\fiverm .})

\plot -3.00000 -3.00000 -1.00000 -1.5 /

\plot -1.00000 -1.5 -2.00000 -3.00000 /

\plot -1.00000 -1.5 0.00000 -3.00000 /

\plot 1.00000 -3.00000 -1.00000 -1.5 /

\putrule from -1.00000 -1.5 to -1.00000 0

\plot -1.00000 0 -3.00000 1.5 /

\plot -2.00000 1.5 -1.00000 0 /

\plot -1.00000 0 0.00000 1.5 /

\plot -1.00000 0 1.00000 1.5 /

\put {$\bullet$} at -3.00000 -3.00000

\put {$\bullet$} at -2.00000 -3.00000

\put {$\bullet$} at 0.00000 -3.00000

\put {$\bullet$} at 1.00000 -3.00000

\put {$\bullet$} at -1.00000 -1.5

\put {$\bullet$} at -1.00000 0

\put {$\bullet$} at -3.00000 1.5

\put {$\bullet$} at -2.00000 1.5

\put {$\bullet$} at 0.00000 1.5

\put {$\bullet$} at 1.00000 1.5

\put{$\cdots$} at -1 -3

\put{$\cdots$} at -1 1.5

\put {$Type \ 3$} at -0.84772 -3.84301
\endpicture
\end{minipage}

\bigskip

We will denote by ${\cal BPK}$ the subvariety of
pseudocomplemented Kleene algebras generated by the subdirectly
irreducible algebras whose dual space is given by one of the above
types and we shall
call algebras in $\mathcal{ BPK}$  {\it bundle} pseudocomplemented
Kleene algebras.

\bigskip

Note that there exist some subdirectly irreducible
pseudocomplemented Kleene algebras whose dual spaces do not have any of the
above types, for example an algebra whose dual space is

$$
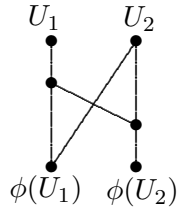
\beginpicture
\setcoordinatesystem units <0.55556cm,0.55556cm> \setlinear
\setshadesymbol <z,z,z,z> ({\fiverm .})

\plot -5 -3 -5 0 /

\plot -5 -3 -3 0 /

\plot -3 -3 -3 0 /

\plot -3 -2 -5 -1 /

\put {$\bullet$} at -5.00000 -3.00000

\put {$\bullet$} at -3.00000 -3.00000

\put {$\bullet$} at -3.00000 -2.00000

\put {$\bullet$} at -5.00000 -1.00000

\put {$\bullet$} at -5.00000 0.00000

\put {$\bullet$} at -3.00000 0.00000

\put {$U_1$} at -5.17677 0.56515 \put {$U_2$} at -2.91617 0.54254
\put {$\phi(U_1)$} at -5.13155 -3.54914 \put {$\phi(U_2)$} at
-2.84835 -3.59435
\endpicture$$

If ${\bf L}_1$ and ${\bf L}_2$ are bounded distributive lattices, let us define
the {\it sum} ${\bf L}_1 \oplus {\bf L}_2$ as the usual {\it ordinal sum} where $1^{{\bf L}_1} = 0^{{\bf L}_2}$.

Our next goal is to show that if ${\bf L}$ is a subdirectly
irreducible $pk$-algebra whose dual space is given by types 1, 2 or 3 then ${\bf L}= {\bf B} \oplus {\bf C} \oplus {\bf B}$, where ${\bf B}$ is a
Boolean  algebra and ${\bf C}$ is a chain. In order to prove this, we will use the term $C(x)$ defined in the previous section along with the term $T(x)$ given by:

$$T(x)=C(x) \wedge C(x)'.$$

Next we show that the variety ${\cal BPK}$ is characterized by the identity $x^\ast \leq C(y) \vee T(x)^\ast$. In order to prove  this  we need the following results.

\begin{Lemma} \label{comp conex con min no incluidos en decreciente 1}
Let ${\bf L}$ be a subdirectly irreducible $pk$-algebra and $V \in \mathbb{D}(\mathbb{X}({\bf L}))$. If $S_i$ is a connected component  of $\mathbb{X}({\bf L})$ such that $V \cap S_i \ne \emptyset $ and $Min(S_i) \not \subseteq (V]$ then there exists a minimal element $M \in \mathbb{X}({\bf L})$ such that $M \in V^\ast$ and $M \not \in T(V)^\ast$.
\end{Lemma}
\begin{Proof}
Assume $V$ and $S_i$ as in the statement of the theorem and consider the following partition of $Min(S_i)$ given by the  sets:

$$A=\{\phi(W) \in S_i: \phi(W) \in (V]\ \} \mbox{ \ \ and\ \ } B=\{\phi(W) \in S_i: \phi(W) \not \in (V] \ \}.$$

Since $[B)$ is not a connected component, there exists $\phi(U) \in A$ and $\phi(R) \in B$ such that $\phi(R) \subseteq U$. Moreover, from $ U \not \in V$, it follows that there exists $T \in V$ such that $\phi(U) \subseteq T$ and $\phi(T) \not \in V$. This is shown in the following picture.

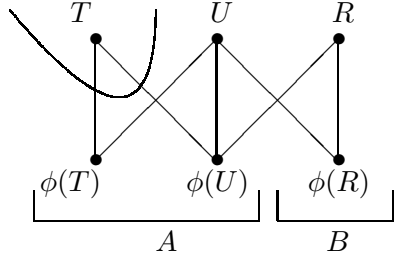
\begin{center}

 \setlength{\unitlength}{0.8cm}
\begin{picture}(6,2.5)
\put(0, 0){\line(0, 1){2}} \put(0, 2){\line(1, -1){2}} \put(2,
0){\line(1, 1){2}} \put(2, 2){\line(1, -1){2}} \put(4, 0){\line(0,
1){2}} \put(2, 0){\line(0, 1){2}} \put(0, 0){\line(1, 1){2}} \put(-.1, -.1){$\bullet$} \put(1.9, -.1){$\bullet$} \put(3.9,
-.1){$\bullet$}
\put(-.1, 1.9){$\bullet$} \put(1.9, 1.9){$\bullet$} \put(3.9,
1.9){$\bullet$}
\put(-.9,-.5){$\phi(T)$} \put(1.5,-.5){$\phi(U)$}
\put(3.5,-.5){$\phi(R)$}
\put(-.4,2.3){$T$} \put(1.9,2.3){$U$} \put(3.9,2.3){$R$}

\put(-1, -1){\line(1, 0){3.7}}
\put(3, -1){\line(1, 0){1.9}}
\put(-1, -1){\line(0, 1){0.5}}
\put(2.7, -1){\line(0, 1){0.5}}
\put(3, -1){\line(0, 1){0.5}}
\put(4.9, -1){\line(0, 1){0.5}}

\put(1,-1.5){$A$}
\put(3.8,-1.5){$B$}

\qbezier(-1.4,2.5)(1,-.4)(1,2.5)

\end{picture}

\end{center}

\

\

Then, using the Lemmas \ref{V inters V'},  and the fact that $Max  \subseteq C(V)$, we have that $\phi(R) \not \in T(V)^\ast$. Moreover, $\phi(R) \in V^\ast$ which completes the proof.

\end{Proof}

\begin{Lemma} \label{comp conex con min no incluidos en decreciente 2}
Let ${\bf L}$ be a subdirectly irreducible $pk$-algebra and $V \in \mathbb{D}(\mathbb{X}({\bf L}))$ such that $V \cap V^\prime \ne \emptyset$. If there exists a connected component $S_{i_0}$ such that $(V \cap V^\prime) \cap S_{i_0} = \emptyset$, then there exists a connected component $S_{j_0}$ such that $(V \cap V^\prime) \cap S_{j_0} \ne \emptyset$ and $Min(S_{j_0}) \not \subseteq (V \cap V^\prime]$.

\end{Lemma}

\begin{Proof}
Suppose that there exists $S_{i_0}$ such that $S_{i_0} \cap (V \cap V^\prime) = \emptyset$ and suppose, contrary to our claim, that $Min(S_j) \subseteq (V \cap V^\prime]$ for every connected component such that $S_j \cap (V \cap V^\prime) \ne \emptyset$.
For each connected component $S_j$ we have that $S_j \cap (V \cap V^\prime) = \emptyset$ or $S_j \cap (V \cap V^\prime) \ne \emptyset$. In the former case we have that $S_j \subseteq (V \cap V^\prime)^\ast$ and so $S_j= \phi(S_j) \subseteq C(V)$. Hence, $S_j \subseteq T(V)^\ast.$ On the contrary, in the latter case, we claim that $S_j \cap T(V)^\ast = \emptyset$. Indeed, under our assumptions, $Min(S_j) \subseteq (V \cap V^\prime]$ and so $Min(S_j) \cap C(V) = \emptyset$. Moreover, from $Max \subseteq C(V)$, $S_j \cap T(V)^\ast = \emptyset$.   These cases show that $T(V)^\ast$ is a union of connected components and $T(V)^\ast \ne \emptyset$. This implies that $T(V)^\ast$ and $(T(V)^\ast)^c$ are Boolean elements, which contradicts the fact that ${\bf L}$ is an indecomposable algebra.

\end{Proof}

Note that it is possible for a subdirectly irreducible $pk$-algebra  to have elements that are maximal and minimal at the same time. This is shown in the following example.

\begin{Example} \label{ejemplo de que no existe el menor elemento denso}
Consider $ X=(\mathbb{N} \times \{0,1\}) \cup \{\infty\}$ with the one-point compactification of the countable discrete space $\mathbb{N} \times \{0,1\}$. The following picture gives an order $\leq$ on $X$. Here $n$ and $n^\prime$ are defined in the same way as in Example \ref{si infinita pero no simple}.

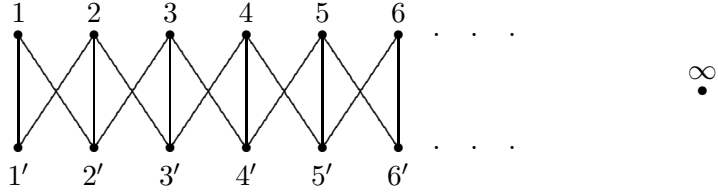
$$
\beginpicture
\setcoordinatesystem units <0.50000cm,0.50000cm>
\setlinear
\setshadesymbol <z,z,z,z> ({\fiverm .})
\putrule from -9.00000 0.00000 to -9.00000 -3.00000
\putrule from -7.00000 0.00000 to -7.00000 -3.00000
\putrule from -5.00000 0.00000 to -5.00000 -3.00000
\putrule from -3.00000 0.00000 to -3.00000 -3.00000
\putrule from -1.00000 0.00000 to -1.00000 -3.00000
\putrule from 1.00000 0.00000 to 1.00000 -3.00000
\plot -9.00000 -3.00000 -7.00000 0.00000 /
\plot -9.00000 0.00000 -7.00000 -3.00000 /
\plot -7.00000 -3.00000 -5.00000 0.00000 /
\plot -7.00000 0.00000 -5.00000 -3.00000 /
\plot -5.00000 0.00000 -3.00000 -3.00000 /
\plot -5.00000 -3.00000 -3.00000 0.00000 /
\plot -3.00000 0.00000 -1.00000 -3.00000 /
\plot -3.00000 -3.00000 -1.00000 0.00000 /
\plot -1.00000 0.00000 1.00000 -3.00000 /
\plot -1.00000 -3.00000 1.00000 0.00000 /
\put {$\scriptstyle\bullet$} at -9.00000 -3.00000
\put {$\scriptstyle\bullet$} at -9.00000 0.00000
\put {$\scriptstyle\bullet$} at -7.00000 -3.00000
\put {$\scriptstyle\bullet$} at -7.00000 0.00000
\put {$\scriptstyle\bullet$} at -5.00000 -3.00000
\put {$\scriptstyle\bullet$} at -3.00000 -3.00000
\put {$\scriptstyle\bullet$} at -1.00000 -3.00000
\put {$\scriptstyle\bullet$} at 1.00000 -3.00000
\put {$\scriptstyle\bullet$} at -5.00000 0.00000
\put {$\scriptstyle\bullet$} at -3.00000 0.00000
\put {$\scriptstyle\bullet$} at -1.00000 0.00000
\put {$\scriptstyle\bullet$} at 1.00000 0.00000
\put {$\scriptstyle\bullet$} at 9.00000 -1.50000
\put {${\bf .}$} at 2.00000 -3.00000
\put {${\bf .}$} at 3.00000 -3.00000
\put {${\bf .}$} at 4.00000 -3.00000
\put {${\bf .}$} at 2.00000 0.00000
\put {${\bf .}$} at 3.00000 0.00000
\put {${\bf .}$} at 4.00000 0.00000
\put {$1'$} at -8.99550 -3.7
\put {$2'$} at -7.01649 -3.7
\put {$3'$} at -5.00750 -3.7
\put {$1$} at -8.98051 0.6
\put {$2$} at -7.00150 0.6
\put {$3$} at -4.99250 0.6
\put {$4$} at -3 0.6
\put {$5$} at -1 0.6
\put {$6$} at 1 0.6

\put {$4'$} at -3 -3.7
\put {$5'$} at -1 -3.7
\put {$6'$} at 1 -3.7

\put {$\infty$} at 9 -1
\endpicture
$$

It is easy to check that $\langle X, \tau, \leq \rangle$ is a $pk$-space which has an element that is maximal and minimal at the same time.

\end{Example}

Nevertheless, this is not possible for a subdirectly irreducible $pk$-algebra that satisfies $x^\ast \leq C(y) \vee T(x)^\ast$.

\begin{Proposition}

Let ${\bf L}$ be a $pk$-algebra that satisfies $x^\ast \leq C(y) \vee T(x)^\ast$ then $\mathbb{X}({\bf L})$ is unitary or $Max(\mathbb{X}({\bf L})) \cap Min(\mathbb{X}({\bf L})) = \emptyset$.

\end{Proposition}

\begin{Proof}
Suppose that $\mathbb{X}({\bf L})$ is not unitary, $U \in \mathbb{X}({\bf L})$ and $U \in Max \cap Min$. Then, $\phi(U)=U$ and since $\mathbb{X}({\bf L})$ is not unitary, there exists a maximal element $U_1$ of $\mathbb{X}({\bf L})$ such that $U_1 \ne U$. Let us consider the clopen increasing $V_{U_1}$ such that $U_1 \in V_{U_1}$ and $U \not \in V_{U_1}$. Since $V_{U_1}$ satisfies the  hypotheses of Lemma \ref{comp conex con min no incluidos en decreciente 2}, there exists a connected component $S_i$ such that $Min(S_i) \not \subseteq (V_{U_1} \cap V_{U_1}^\prime]$.
Applying Lemma \ref{comp conex con min no incluidos en decreciente 1} to $V_{U_1} \cap V_{U_1}^\prime$ we have that there exists a minimal element $\phi(R)$ such that $\phi(R) \in (V_{U_1} \cap V_{U_1}^\prime)^\ast$ and $\phi(R) \not \in T(V_{U_1} \cap V_{U_1}^\prime)^\ast$. Again, since $ R \not \subseteq \phi(R)$, there exists a clopen increasing $V_R$ such that $R \in V_R$ and $\phi(R) \not \in V_R$. Finally, by Lemma \ref{V inters V'}, we obtain that $\phi(R) \not \in C(V_R)$. From this, $\phi(R) \in (V_{U_1} \cap V_{U_1}^\prime)^\ast$ but $\phi(R) \not \in T(V_{U_1} \cap V_{U_1}^\prime)^\ast \vee C(V_R)$, which contradicts our assumption.

\end{Proof}

Clearly, the previous lemma and the fact that $Max$ is a closed set guarantee the existence of the least dense element in a subdirectly irreducible $pk$-algebra satisfying $x^\ast \leq C(y) \vee T(x)^\ast$, namely, the element associated with the clopen increasing set $Max$. Observe that Example \ref{ejemplo de que no existe el menor elemento denso} shows that there exist subdirectly irreducible $pk$-algebras that do not have a least dense element. It is enough to consider the sequence of dense elements: $V_1=Max \cup \{2^\prime, 3^\prime, \ldots \}$, $V_2=Max \cup \{3^\prime, 4^\prime, \ldots \}$, $V_3=Max \cup \{4^\prime, 5^\prime, \ldots \}$, etc.

Here, and subsequently, $d$ denotes the least dense element whenever it exists.

\begin{Lemma} \label{T}
Let ${\bf L}$ be a subdirectly irreducible $pk$-algebra that satisfies $x^* \leq C(y) \vee T(x)^*$. Then for every
$a  \in L$, $a \ne 0$, we have that $a^* \leq d$.
\end{Lemma}

\begin{Proof}
For $a = 1$ the result is trivial. Let $a \in L$, $a \ne 0, 1$. Let
$V$ be the clopen increasing set that represents $a$. Suppose, contrary to our claim, that there exists a prime filter $Q$ in $V^* = (V]^c$, $Q \not
\in Max$. We will prove that there exists $R \in
\mathbb{X}({\bf L}) \setminus Max$ such that $R \in
V^*$ but $R \not \in T(V)^*$.

Let $\mathbb{X}({\bf L}) = \bigcup_{i \in I} S_i$ be the decomposition of $\mathbb{X}({\bf L})$ in connected components. The proof will be divided in two cases.

\begin{itemize}
\item[(a)] Suppose that $Min \subseteq (V]$.
Since $Q \in V^*$ and $Q \not \in Max$, then
$Q \in Body({\bf L})$. Therefore, $Q=P$ or $Q=\phi(P)$ and thus
$P \in V^*$. Let us prove now that $P \not \in T(V)^*$.

Let $U \in Max$ such that $P \subseteq U$. From $P \in V^\ast$, we have that $U \not \in V$ and $\phi(U) \not \in V$. Since $\phi(U)$ is minimal and $Min \subseteq (V]$, there exists $R \in V$ such that $\phi(U) \subseteq R$. Moreover, from $U \not \in V$, $\phi(R) \not \in V$. This is shown in the following picture.

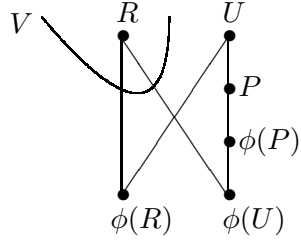
\begin{center}

 \setlength{\unitlength}{0.7cm}
\begin{picture}(3,3)
\hspace{-1cm}\put(0.1, -.9){\line(0, 1){3}}
\put(2.1, -.9){\line(0, 1){3}}
\put(2.1, -.9){\line(-2, 3){2}}
 \put(0.1, -.9){\line(2, 3){2}}
\put(0, -1){$\bullet$}
\put(2, -1){$\bullet$}
\put(0,2){$\bullet$}
\put(2, 0){$\bullet$} \put(2, 1){$\bullet$}
\put(2,2){$\bullet$}
\put(-0.1,-1.5){$\phi(R)$} \put(2,-1.5){$\phi(U)$}
\put(2.3,0){$\phi(P)$}
\put(2.3,1){$P$}
\put(0,2.4){$R$}
\put(2,2.4){$U$}
\put(-2,2.2){$V$}

\qbezier(-1.4,2.5)(1,-.4)(1,2.5)

\end{picture}

\end{center}

\

\

From Lemma \ref{condiciones de V y V prima} it is immediate that $U \in T(V)$ and so, from $P \subseteq U$, we have that $P \not \in T(V)^\ast$.

\item[(b)] Suppose that $Min \not \subseteq (V]$. Again, the proof will be divided in two cases.

\begin{enumerate}

\item[$(i)$] Suppose that $V \cap  S_i \ne \emptyset$ for each connected component $S_i$ of $\mathbb{X}({\bf L})$. Under the above assumptions, there exists $S_i$  such that $S_i \cap V \ne \emptyset$ and $Min(S_i) \not \subseteq (V]$. Then, by Lemma \ref{comp conex con min no incluidos en decreciente 1}, there exists a minimal element $\phi(R)$ such that $\phi(R) \in V^\ast$ and  $\phi(R) \not \in T(V)^\ast$.

\item[$(ii)$] Suppose that there exists a connected component $S_j$ such that $V \cap  S_j = \emptyset$. Applying  Lemma \ref{comp conex con min no incluidos en decreciente 2} to the clopen increasing $V$, there exists a connected component $S_i$ such that $Min(S_i) \not \subseteq (V \cap V^\prime]$. It is enough to apply  Lemma \ref{comp conex con min no incluidos en decreciente 1} to the clopen increasing $W=V \cap V^\prime$ to find a minimal element $\phi(R)$ such that $\phi(R) \in W^\ast$ and $\phi(R) \not \in T(W)^\ast$.
\end{enumerate}

\end{itemize}

These two cases show that if there exists  $V \in
\mathbb{D}(\mathbb{X}({\bf L}))$, $V \ne \emptyset, V \ne \mathbb{X}({\bf L})$ such that $V^*  \not
\subseteq Max$, then  there exists $R \in
\mathbb{X}({\bf L}) \setminus Max$ such that $R \in
V_1^*$ but $R \not \in T(V_1)^*$ for some clopen increasing $V_1$.
In addition if we take $W = \eta(d) =  Max$, we have that $C(W)=Max$.
Therefore, taking $W=\eta(d)=Max$, we obtain $V_1^* \not \subseteq C(W) \vee T(V_1)^*$.

We have thus proved that if ${\bf L}$ satisfies $x^* \leq C(y) \vee
T(x)^*$ then for every $a \in L$, $a \ne 0$, we have that  $a^* \leq
d$.
\end{Proof}

\begin{Lemma}
Let ${\bf L}$ be a subdirectly irreducible $pk$-algebra that satisfies $x^*  \leq C(y) \vee T(x)^*$. Then $ L = [0,
d] \cup [d, d'] \cup [d', 1]$.
\end{Lemma}

\begin{Proof}
If $d=1$ the proof es trivial. Suppose that $d \ne 1$. Observe
that for every $a \in L$,  either $a \leq d$ or $a \geq d$. Indeed,
if we suppose that $a \not \geq d$, then $a^* \ne 0$. But $a \leq
a^{**} $ and by Lemma \ref{T}, $a^{**} \leq d$, so $a \leq d$.

Now we apply this to both $a$ and $a'$. If $a \leq d$ then $a \in
[0,d]$. If $a \geq d$ and $a' \leq d$ then $a \geq d'$ and so $a \in
[d',1]$,  and if $a \geq d$ and $a' \geq d$ then $a \geq d$ and $a
\leq d'$, and so $a \in [d,d']$. So $a \in [0, d] \cup [d, d'] \cup
[d', 1]$.\end{Proof}

\medskip

\begin{Proposition}
Let ${\bf L}$ be a subdirectly irreducible $pk$-algebra that satisfies $x^* \leq C(y) \vee T(x)^*$. Then  $[d, d']$
has at most three elements.
\end{Proposition}
\begin{Proof}
From the duality, $\eta(d) = Max$ and
$\eta(d')=\mathbb{X}({\bf L}) \setminus Min$, thus
if  $a \in [d,d']$ we have that $Max \subseteq
\eta(a) \subseteq \mathbb{X}({\bf L}) \setminus
Min$.

On the other hand, by Theorem \ref{body de s.i.},
$|\mbox{Body}({\bf L})| \leq 2$ and $\phi(P)$ and $ P$ are
comparable elements, for each filter $P$, then either $\eta(a)=
Max$ or $\eta(a) = Max \cup \{P\}$ or $\eta(a)= Max \cup \{P,
\phi(P)\} = \mathbb{X}({\bf L}) \setminus Min$.
Consequently, $[d, d']$ is a chain with  1, 2 or 3 elements.
\end{Proof}

\

It is well known that the class of pseudocomplemented algebras
satisfies  Glivenko's Theorem (see \cite{Bal-Dwi}), that is,
{\cal Rg}$({\bf L}) \cong {\bf L}/{\cal D}({\bf L})$, where {\cal Rg}$({\bf L})$ is the
Boolean algebra of regular elements and ${\cal D}({\bf L})$  is the filter
of dense elements of ${\bf L}$. We showed that in a
subdirectly irreducible $pk$-algebra ${\bf L}$ satisfying $x^\ast \leq C(y) \vee T(x)^\ast$ there
exists a least dense element $d$, that corresponds (under the duality)
to the clopen increasing set $Max$. Consequently, we have that ${\cal
D}({\bf L}) =[d, 1]$. Therefore, by the above, we obtain that
$$[0, d] \cong {\bf L}/ {\cal D}({\bf L}),$$ which is a Boolean algebra.

Furthermore, taking into account that $'$ is an involutive
anti-isomorphism, it is clear that $[d',1]$ has also an underlying
structure of Boolean lattice.
So we have the following theorem.

\medskip

\begin{Theorem} \label{forma de si en bpk}
If ${\bf L}$ is a subdirectly irreducible $pk$-algebra that satisfies $x^*  \leq C(y) \vee T(x)^*$, then ${\bf L}$ has an
underlying structure of distributive lattice given by $${\bf B} \oplus {\bf C}_i
\oplus {\bf B},$$ where ${\bf C}_i$ is an $i$-element chain, for $1 \leq i \leq
3$, and ${\bf B}$ is a Boolean lattice.
\end{Theorem}

\begin{Corollary} \label{min en max}
Let ${\bf L}$ be a subdirectly irreducible $pk$-algebra that satisfies $x^*  \leq C(y) \vee T(x)^*$. If $P \in
\mathbb{X}({\bf L}) \setminus Max(\mathbb{X}({\bf L}))$ then $P \subseteq U$ for
every $U \in Max$. Moreover, $\mathbb{X}({\bf L})$ has Type 1, Type 2 or Type 3.
\end{Corollary}
\medskip

The above results together with the next proposition provide an
equational basis for the variety $\mathcal{ BPK}$ of bundle pseudocomplemented Kleene algebras.

\begin{Proposition} \label{formas implican ecuacion general}
If ${\bf L}$ is a subdirectly irreducible $pk$-algebra such that $\mathbb{X}({\bf L})$ has Type 1, 2 or 3, then ${\bf L}$
satisfies $x^* \leq C(y) \vee T(x)^*$.

\end{Proposition}

\begin{Proof}
Let $a, b \in L$, $V=\eta(a)$ and $W= \eta(b)$ the clopen increasing
sets that represent $a$ and $b$ respectively. Suppose that $V \ne
\emptyset$. Then, since $\mathbb{X}({\bf L})$ is given by Type 1, Type 2 or
Type 3, we obtain $(V]^c \subseteq Max(\mathbb{X}({\bf L}))$.

Furthermore, since $C(b)$ is a dense element, we have that
$Max \subseteq C(W)$. From this, $V^*
\subseteq C(W)$.

Finally, if $V= \emptyset$, by Lemma \ref{T} it follows that $V^* =
\mathbb{X}({\bf L})$ and $T(V)^\ast=\mathbb{X}({\bf L})$ which establishes $V^*
\subseteq C(W) \cup T(V)^*$.
\end{Proof}

\


To close this section we summarize here the most important results obtained in the last two sections.

\begin{itemize}
\item Equational bases for the varieties ${\cal PK}_0$, ${\cal PK}_1$ and ${\cal BPK}$ relative to the variety ${\cal PK}$:

\begin{itemize}
\item ${\bf L} \in \mathcal{ PK}_0 \Longleftrightarrow {\bf L} \models C(x)^\prime \leq C(x)$
\item ${\bf L} \in \mathcal{ PK}_1 \Longleftrightarrow {\bf L} \models C(x) \wedge C(x)^\prime \leq C(y)$
\item ${\bf L} \in \mathcal{ BPK} \Longleftrightarrow {\bf L} \models x^\ast \leq  C(y) \vee T(x)^\ast $
\end{itemize}

\item Description of the subdirectly irreducible members of ${\cal BPK}$:
\begin{itemize}
\item ${\bf L}$ is s.i. in ${\cal BPK}$  $\Longleftrightarrow$  $\mathbb{X}({\bf L})$ has Type 1, 2 or 3 \ $\Longleftrightarrow$ ${\bf L} \cong {\bf B} \oplus {\bf C}_i \oplus {\bf B}$, where ${\bf B}$ is a boolean algebra and ${\bf C}$ is a chain of at most 3 elements.
\item ${\bf L}$ s.i. in ${\cal BPK}, {\bf L}$ non trivial $\Longrightarrow$ $Min \cap Max = \emptyset$
\end{itemize}



\end{itemize}

\section{Subvarieties of ${\cal BPK}$}

The objective of this section is to determine the lattice of
subvarieties of $\mathcal{ BPK}$, the subvariety of those
pseudocomplemented Kleene algebras that satisfy $x^* \leq C(y) \vee
T(x)^*$.

 The following results are easy to check.

\begin{Lemma} \label{Bolle n generada equivale a Bk n generada}

\item Let ${\bf L}={\bf B} \oplus {\bf C}_m \oplus {\bf B} $ be a subdirectly irreducible algebra in ${\cal BPK}$, then ${\bf L}$ satisfies the following.
\begin{enumerate}
\item If ${\bf S}$ is a Boolean sublattice of ${\bf B}$, then $S \cup S^\prime$ is a subalgebra of ${\bf L}$.

\item If ${\bf S}$ is a subalgebra of ${\bf L}$ and $S \ne \{0,1\}$, then ${\bf S} \cap {\bf B}$ is a subalgebra of the Boolean algebra ${\bf B}$.
\end{enumerate}

\end{Lemma}

Given an algebra ${\bf L}$ in ${\cal BPK}$ and a subset $G  \subseteq L$, we denote by $G^\prime=\{g^\prime: g \in G\}$ and by $Sg(G)$ the subalgebra generated by $G$.
If $G$ is a finite subset of ${\bf L} = [0,d] \oplus {\bf C}_m \oplus [d^\prime, 1]$ with $m \in \{1,2,3\}$ and we consider the set $\tilde{G}_1= \tilde{G} \cap [0,d]$, where $\tilde{G}= G \cup G^\prime$ by the lemma above it follows that $Sg(G)=S \cup S^\prime$ or $Sg(G)= S \cup S^\prime \cup C_m$, where $S$ is the Boolean lattice generated by $\tilde{G}_1$. Hence, since Boolean algebras are locally finite, it follows that $Sg(G)$  is a subdirectly irreducible finite algebra. Furthermore, we know that the $n$-generated free algebra ${\cal F}_{{\cal BPK}}(n)$ is isomorphic to a subdirect product of  the members of  $\{{\cal F}_{{\cal BPK}}(n)/\theta, \theta \in {\cal M} \}$, where ${\cal M}$  is the set of maximal congruences on ${\cal F}_{{\cal BPK}}(n)$. Since  ${\cal F}_{{\cal BPK}}(n)/\theta$ is a finite algebra and ${\cal M}$ is a finite set, we have that ${\cal F}_{{\cal BPK}}(n)$ is a finite algebra.  This shows that the variety ${\cal BPK}$ is locally finite and so it is generated by its finite subdirectly irreducible members.

Here and subsequently, we will denote a finite subdirectly irreducible algebra ${\bf L}$ in ${\cal BPK}$  by $\mathbb{B}_{(i,m)}= \mathbb{B}^-_i \oplus {\bf C}_m \oplus \mathbb{B}^+_i $, where $ \mathbb{B}^-_i$ is the sublattice  of ${\bf L}$ with $i$ atoms given by $[0,d]$, $ \mathbb{B}^+_i$ is the sublattice of ${\bf L}$ with $i$ atoms given by $[d^\prime,1]$ and ${\bf C}_m$ es the chain $[d, d^\prime]$.  Recall that $ \mathbb{B}^-_i$ and $\mathbb{B}^+_i$ have an underlying structure of Boolean lattice.

Given a class ${\cal K}$ of algebras,  ${\bf Si}({\cal K})$
 consists of precisely one algebra from
each of the isomorphism classes of the subdirectly irreducible algebras. Let ${\bf Si}_{{\bf fin}}({\cal K})$ be the class of finite members of ${\bf Si}({\cal K})$. For ${\cal BPK}$, observe that $ {\bf Si}_{{\bf fin}}({\cal BPK})= \{\mathbb{B}_{(i,m)}: i \in \mathbb{N}_0, m \in \{1,2,3\}\}$.

Since $\mathcal{BPK}$ is congruence-distributive and locally finite
we can apply   J\'onsson's Theorem  (see \cite{Jonsson}) and its
generalization by Davey (see \cite{Davey}) to find the lattice $\Lambda(\mathcal{
BPK})$ of subvarieties of $\mathcal{BPK}$. So, $\Lambda({\cal BPK})$ is a complete distributive lattice isomorphic to  the lattice of the decreasing sets of the order set ${\bf Si}_{{\bf fin}}({\cal BPK})$. Recall that the order defined on ${\bf Si}_{{\bf fin}}({\cal BPK})$ is given by ${\bf A} \preceq {\bf B} \Leftrightarrow {\bf A} \in HS({\bf B})$.

Now we are going to find a simpler way of determining the order relation on
${\bf Si_{{\bf fin}}}(\mathcal{ BPK})$.

\begin{Proposition}
If $\mathbb{B}_{(i,m)}, \mathbb{B}_{(j,n)} \in {\bf Si}_{{\bf fin}}(\mathcal{ BPK})$ then
$$\mathbb{B}_{(i,m)} \preceq \mathbb{B}_{(j,n)} \ \Longleftrightarrow \ \mathbb{B}_{(i,m)} \in
HS(\mathbb{B}_{(j,n)}) \  \Longleftrightarrow \ i \leq j \mbox{ and } m \leq n$$
\end{Proposition}

\begin{Proof}
It is immediate that if $\mathbb{B}_{(i,m)} \in HS(\mathbb{B}_{(j,n)})$ then $i \leq
j$ and $m \leq n$.

For the converse, consider the following cases: If $m=n$ and $i \leq
j$ then $\mathbb{B}_{(i,m)}$ is a subalgebra of $\mathbb{B}_{(j,m)}$, so $\mathbb{B}_{(i,m)}
\in HS(\mathbb{B}_{(j,m)})$. If $m=2$ and $n=3$ and $i \leq j$ then $\mathbb{B}_{(i,m)}$ is
a subalgebra of $\mathbb{B}_{(j,n)}$, so $\mathbb{B}_{(i,m)} \in HS(\mathbb{B}_{(j,n)})$. If
$m=1 <n$ and $i \leq j$, taking into account that $\mathbb{B}_{(i,1)} \in H(\mathbb{B}_{(i,n)})$ (we can identify $[d, d^\prime]$ with the element $d=d^\prime$ in $\mathbb{B}_{(i,1)}$), we obtain $\mathbb{B}_{(i,m)} \in HS(\mathbb{B}_{(j,n)})$.
\end{Proof}
\medskip

The next picture shows the ordered set  ${\bf Si}_{{\bf fin}}(\mathcal{BPK})$. Here, $\mathbb{B}_{(1,0)}$ denotes the chain with two elements and $\mathbb{B}_{(0,0)}$ is the trivial algebra.
\medskip

$$\beginpicture
\setcoordinatesystem units <0.80000cm,0.80000cm>
\setlinear
\setshadesymbol <z,z,z,z> ({\fiverm .})
\plot -7.00000 -6.00000 -5.00000 -7.00000 /
\plot -5.00000 -7.00000 -3.00000 -8.00000 /
\plot -7.00000 -4.00000 -3.00000 -6.00000 /
\plot -7.00000 -2.00000 -3.00000 -4.00000 /
\plot -7.00000 0.00000 -3.00000 -2.00000 /
\putrule from -7.00000 0.00000 to -7.00000 -6.00000
\putrule from -5.00000 -1.00000 to -5.00000 -7.00000
\putrule from -3.00000 -2.00000 to -3.00000 -8.00000
\putrule from -3.00000 -8.00000 to -3.00000 -9.00000
\putrule from -3.00000 -8.00000 to -3.00000 -10.00000

\put {$\bullet$} at -3.00000 -10.00000
\put {$\bullet$} at -3.00000 -8.00000
\put {$\bullet$} at -5.00000 -7.00000
\put {$\bullet$} at -5.00000 -5.00000
\put {$\bullet$} at -3.00000 -6.00000
\put {$\bullet$} at -7.00000 -4.00000
\put {$\bullet$} at -7.00000 -6.00000
\put {$\bullet$} at -7.00000 -2.00000
\put {$\bullet$} at -5.00000 -3.00000
\put {$\bullet$} at -3.00000 -4.00000
\put {$\bf{.}$} at -7.00000 0.00000
\put {$\bullet$} at -5.00000 -1.00000
\put {$\bullet$} at -3.00000 -2.00000
\put {$\bullet$} at -7.00000 0.00000
\put {${\bf .}$} at -7.00000 1.00000
\put {${\bf .}$} at -7.00 1.5
\put {${\bf .}$} at -5.00000 0.00000
\put {${\bf .}$} at -5.00000 0.00000
\put {${\bf .}$} at -5 0.5
\put {${\bf .}$} at -3.00000 -1.00000
\put {${\bf .}$} at -3.00000 -1.00000
\put {$\bullet$} at -3.00000 -9.00000
\put {${\bf .}$} at -7 0.5
\put {${\bf .}$} at -5 -0.5
\put {${\bf .}$} at -3 -1.5
\put {${\bf .}$} at -3 -0.5
\put {$\mathbb{B}_{(0,0)}$} at -2 -10.32534
\put {$\mathbb{B}_{(1,0)}$} at -2 -9.32534
\put {$\mathbb{B}_{(1,1)}$} at -2 -7.94603
\put {$\mathbb{B}_{(2,1)}$} at -2 -6.20690
\put {$\mathbb{B}_{(3,1)}$} at -2 -4.25787
\put {$\mathbb{B}_{(4,1)}$} at -2 -2.30885
\put {$\mathbb{B}_{(1,2)}$} at -4.07796 -6.92654
\put {$\mathbb{B}_{(2,2)}$} at -4.01799 -4.91754
\put {$\mathbb{B}_{(3,2)}$} at -4.13793 -2.78861
\put {$\mathbb{B}_{(4,2)}$} at -3.92804 -0.95952
\put {$\mathbb{B}_{(1,3)}$} at -6.14693 -5.69715
\put {$\mathbb{B}_{(2,3)}$} at -6.08696 -3.71814
\put {$\mathbb{B}_{(3,3)}$} at -6.05697 -1.79910
\put {$\mathbb{B}_{(4,3)}$} at -6.11694 0.26987
\endpicture$$

Let ${\cal V}(\mathbb{B}_{(i,m)})$ denote the variety generated by
$\mathbb{B}_{(i,m)}$ and  let us denote by ${\cal BPK}_0$ and ${\cal BPK}_1$ the varieties
${\cal BPK}_0= {\cal V}\left(\{\mathbb{B}_{(i,1)}: i \geq 1\}\right)$ and ${\cal
BPK}_1={\cal V}\left(\{\mathbb{B}_{(i,1)}:i \geq 1\}  \cup
\{\mathbb{B}_{(i,2)}: i \geq 1\}\right)$. Clearly, ${\cal BPK}_0 \subseteq
{\cal BPK}_1 \subseteq \mathcal{BPK}$.

\medskip

We know by J\'onsson's Theorem  that ${\cal V}(\mathbb{B}_{(i,m)})$ is
join-irreducible for every pair $(i,m)$, $1 \leq m \leq 3$.

Let us now find the set ${\cal J}(\Lambda(\mathcal{ BPK}))$ of all
join-irreducible elements of the lattice $\Lambda(\mathcal{ BPK})$ of
all subvarieties of $\mathcal{ BPK}$.
\medskip

\begin{Proposition}
${\cal BPK}_0$, ${\cal BPK}_1$ and $\mathcal{ BPK}$ are join-irreducible elements of $\Lambda(\mathcal{ BPK})$.
\end{Proposition}

\begin{Proof}
Suppose that ${\cal BPK}_0={\cal V}_1 \vee {\cal V}_2$. Then if ${\bf L}$
is a finite subdirectly irreducible algebra of ${\cal V}_1$ or
${\cal V}_2$, we have that ${\bf L}=\mathbb{B}_{(i,1)}$ for some $i$. Let $I_1=\{i:
\mathbb{B}_{(i,1)} \in {\cal V}_1\}$ and $I_2=\{i: \mathbb{B}_{(i,1)} \in {\cal
V}_2\}$. If both $I_1$ and $I_2$ are finite then there exists $i_0$
such that $\mathbb{B}_{(i_0,1)} \not \in {\cal V}_1$ and $\mathbb{B}_{(i_0,1)} \not \in
{\cal V}_2$, which is impossible. Therefore, $I_1$ is infinite or
$I_2$ is infinite and thus ${\cal V}_1={\cal BPK}_0$ or ${\cal
V}_2={\cal BPK}_0$.

A similar argument shows that ${\cal BPK}_1$ and $\mathcal{ BPK}$ are
join-irreducible.
\end{Proof}

\begin{Theorem}
${\cal BPK}_0$, ${\cal BPK}_1$ and $\mathcal{ BPK}$ are the unique
join-irreducible varieties that are not finitely generated.
\end{Theorem}

\begin{Proof}
Let ${\cal V}$ be a join-irreducible non finitely generated variety
of $\mathcal{ BPK}$. Consider the following sets
$$I_1=\{i: \mathbb{B}_{(i,1)} \in {\cal V}\}, \ \ I_2=\{i: \mathbb{B}_{(i,2)} \in {\cal V}\} \ \ I_3=\{i: \mathbb{B}_{(i,3)} \in {\cal V}\}$$
and consider the varieties associated with them
$${\cal V}_1= {\cal V}(\{\mathbb{B}_{(i,1)}, i \in I_1 \}), \ \ {\cal V}_2={\cal V}(\{\mathbb{B}_{(i,2)}, i \in I_2 \}), \ \ {\cal V}_3={\cal V}(\{\mathbb{B}_{(i,3)}, i \in I_3 \}).$$
Thus,  ${\cal V}={\cal V}_1 \vee {\cal V}_2 \vee {\cal V}_3$. If
$I_1, I_2$ and $I_3$ are finite then ${\cal V}$ is finitely
generated, which is a contradiction. Hence, some of them is
infinite. If $I_3$ is infinite then ${\cal V}_3=\mathcal{ BPK}$, so ${\cal V}=\mathcal{ BPK}$. If $I_2$ is infinite and $I_3$ is
finite then ${\cal V}_1 \vee {\cal V}_2= {\cal BPK}_1$ and if $i_3$
is the biggest element of $I_3$, we have that ${\cal
V}_3={\cal V}(\mathbb{B}_{(i_3,3)})$. From this, we obtain  ${\cal V}={\cal BPK}_1
\vee V(\mathbb{B}_{(i_3,3)})$ and since ${\cal V}$ is not finitely generated
and it is join-irreducible, ${\cal V}={\cal BPK}_1$. We now suppose
that $I_1$ is infinite but $I_2$ and $I_3$ are not and let $i_2$ and
$i_3$ be the biggest elements of $I_2$ and $I_3$ respectively. Then
${\cal V}_1={\cal BPK}_0$, ${\cal V}_2={\cal V}(\mathbb{B}_{(i_2,2)})$ and  ${\cal
V}_3={\cal V}(\mathbb{B}_{(i_3,3)})$. It follows that ${\cal V}={\cal BPK}_0 \vee
{\cal V}(\mathbb{B}_{(i_2,2)}) \vee {\cal V}(\mathbb{B}_{(i_3,3)})$, and since ${\cal V}$ is
join-irreducible, we obtain that ${\cal V}={\cal BPK}_0$.

\end{Proof}

\begin{Remark}
We have shown that the join-irreducible elements in $\Lambda({\cal BPK})$ are the varieties  $\mathcal{ BPK}$, ${\cal BPK}_1$, ${\cal BPK}_0$, ${\cal V}(\mathbb{B}_{(i_1,1)})$,
${\cal V}(\mathbb{B}_{(i_2,2)})$ or ${\cal V}(\mathbb{B}_{(i_3,3)})$.
Moreover, it is easy to see that if ${\cal V}$ is an element of $\Lambda({\cal BPK})$ that is not join-irreducible, ${\cal V}$ has one of the following forms: ${\cal V}(\mathbb{B}_{(i_1,1)}) \vee {\cal V}(\mathbb{B}_{(i_2,2)}) \vee {\cal V}(\mathbb{B}_{(i_3,3)})$,
 ${\cal V}(\mathbb{B}_{(i_1,1)}) \vee {\cal V}(\mathbb{B}_{(i_2,2)})$,
${\cal BPK}_0 \vee {\cal V}(\mathbb{B}_{(i_2,2)}) \vee {\cal V}(\mathbb{B}_{(i_3,3)})$,
 ${\cal BPK}_1 \vee {\cal V}(\mathbb{B}_{(i_3,3)})$,
  ${\cal BPK}_0 \vee {\cal V}(\mathbb{B}_{(i_2,2)})$, ${\cal BPK}_0 \vee {\cal V}(\mathbb{B}_{(i_3,3)})$, in all cases $i_1 >i_2 > i_3$.

\end{Remark}

\section{Equational Bases}

In this section we will find equational bases for each subvariety of
$\mathcal{ BPK}$ obtained in the previous section.  It is immediate, from  Theorem \ref{identidades para PK0 y PK1} that the variety ${\cal BPK}_0$ is the subvariety of ${\cal PK}$ characterized by the identities $x^\ast \leq C(y) \vee T(x)^\ast$ and $C(x)^\prime \leq C(x)$ and the variety ${\cal BPK}_1$ is the subvariety of ${\cal PK}$ characterized by $x^\ast \leq C(y) \vee T(x)^\ast$ and $C(x)^\prime \wedge C(x) \leq C(y)$. Next, we will show that it is possible to characterize the variety ${\cal BPK}_1$ whit only one identity.

\begin{Theorem}
Let ${\bf L}$ be a subdirectly irreducible algebra in ${\cal BPK}$ that
satisfies
$$T(x)^* \vee C(y)=C(x) \vee T(y)^*.$$ Then $\mathbb{X}({\bf L})$ has Type 1 or 2, that is, ${\bf L} \in {\cal BPK}_1$.
\end{Theorem}

\begin{Proof}
Let us prove first  that if ${\bf L} $ satisfies  $T(x)^\ast \vee C(y) =C(x) \vee T(y)^\ast$ for all $x, y \in  L$ then it also satisfies $x^\ast \leq C(y) \vee T(x)^\ast$. Indeed, since $x \wedge x^\prime \leq x$, we have that $x^\ast \leq (x \wedge x^\prime)^\ast \leq C(x)$ and so $x^\ast \leq C(x) \vee T(y)^\ast = C(y) \vee T(x)^\ast$.

Finally let us prove  that if $P \in Body({\bf L})$ then $P = \phi(P)$. By way of contradiction, suppose  that $P \ne \phi(P)$.
We show that ${\bf L}$ does not satisfy the identity
$T(x)^* \vee C(y)=C(x) \vee T(y)^*$. Without loss of generality we
can assume that $\phi(P) \subseteq P$. Let us consider the sets $V = Max$
and $W = Max \cup \{P\}$. Then we have that $T(V)^*
\vee C(W)=Max \cup \{P\}$ and $C(V) \vee T(W)^* =
Max$, which is a contradiction.
\end{Proof}

\begin{Remark}
If ${\bf L}$ is a subdirectly irreducible $pk$-algebra such that $\mathbb{X}({\bf L})$ has Type 1 or 2, then
$C(V)=Max$ for every
 $V \ne \mathbb{X}({\bf L}), V \ne\emptyset$, that is, $C(x)$ is the least dense element of ${\bf L}$ for all $x \in L$, $x \ne 0,1$.
\end{Remark}

\begin{Proposition}
Let ${\bf L}$ be a subdirectly irreducible algebra in $\mathcal{ BPK}$ such
that $\mathbb{X}({\bf L})$ has Type 1 or 2. Then ${\bf L}$ satisfies
$T(x)^* \vee C(y)=C(x) \vee T(y)^*$.
\end{Proposition}

\begin{Proof}
Let $V \ne \mathbb{X}({\bf L}), V \ne \emptyset, W \ne \mathbb{X}({\bf L})$ and $W \ne \emptyset$. Then $T(V)^* = T(W)^* =
\emptyset$ and $C(V)= C(W)= Max$. Thus the
identity  $T(x)^* \vee C(y)=C(x) \vee T(y)^*$ holds for all $x \ne 0$, $x \ne
1 $, $y \ne 0$, $y \ne 1$ in $L$. If $V=\emptyset$ or $V= \mathbb{X}({\bf L})$, then
$C(V)=T(V)^* = \mathbb{X}({\bf L})$. From this, ${\bf L}$ satisfies the identity
$T(x)^* \vee C(y)=C(x) \vee T(y)^*$.
\end{Proof}

\medskip

If we write $\gamma(x,y) = T(x)^* \vee C(y)$, we have that ${\cal BPK}$, ${\cal BPK}_1$, and ${\cal BPK}_0$ are characterized within the variety ${\cal PK}$ by the following identities:

\begin{itemize}
\item $\mathcal{ BPK}$: $x^* \wedge \gamma(x,y) \approx x^\ast$

\item ${\cal BPK}_1$: $\gamma(x,y) \approx \gamma(y,x)$

\item ${\cal BPK}_0$: $x^* \wedge \gamma(x,y) \approx x^\ast$ and $C(x)^\prime \wedge C(x) \approx C(x)^\prime$.

\end{itemize}

We will use the following result (see \cite[p.162]{Bal-Dwi}) to show that if a subdirectly irreducible algebra ${\bf L} \in {\cal BPK}$ satisfies certain identity
then its dual space has at most a fixed number  of maximal elements.

\begin{Theorem}
For a $p$-algebra ${\bf L}$ and $n >0$ the following are equivalent:

\begin{enumerate}
  \item[$(1)$] $(x_0 \wedge \ldots \wedge x_{n-1})^* \vee \bigvee_{i<n} (x_0 \wedge \ldots \wedge x_{i-1} \wedge x_i^* \wedge x_{i+1} \wedge \ldots x_{n-1})^* \approx 1$ is an identity
  in ${\bf L}$.

  \item[$(2)$] Every prime filter is contained in at most $n$ maximal filters.
\end{enumerate}
\end{Theorem}

We denote the identity in the last theorem by $\beta(x_0, x_1, \ldots, x_{n-1}) \approx 1$. Since in a subdirectly irreducible algebra ${\bf L} \in {\cal BPK}$, every maximal element contains all non-maximal elements, it follows immediately by the theorem above, that if ${\bf L}$ satisfies $\beta(x_0, x_1, \ldots, x_{n-1}) \approx 1$, its dual space  has at most $n$ maximal elements. Hence, we may derive the following result.

\medskip

\begin{Theorem}
Given a $pk$-algebra ${\bf L}$:

\begin{itemize}
\item ${\bf L} \in \mathcal{ BPK} \Longleftrightarrow {\bf L} \models x^\ast \wedge \gamma(y,x) \approx x^\ast$,

\item
${\bf L} \in {\cal BPK}_1 \Longleftrightarrow {\bf L} \models  \gamma(y,x) \approx \gamma(x,y)$,

\item ${\bf L} \in {\cal BPK}_0 \Longleftrightarrow  {\bf L} \models  x^\ast \wedge \gamma(y,x) \approx x^\ast \mbox{  \ \ and \  \ } {\bf L} \models  C(x) \wedge C(x)' \approx C(x)'$,

\item ${\bf L} \in {\cal V}(\mathbb{B}_{(n,1)}) \Longleftrightarrow  {\bf L} \models  x^\ast \wedge \gamma(y,x) \approx x^\ast, \ {\bf L} \models  C(x) \wedge C(x)' \approx C(x)'  \mbox{ \ \ and \ \ }$  ${\bf L} \models \beta(x_0, x_1, \ldots, x_{n-1}) \approx 1$,

\item ${\bf L} \in {\cal V}(\mathbb{B}_{(n,2)}) \Longleftrightarrow {\bf L} \models \gamma(y,x) \approx \gamma(x,y) \ \ \mbox{and} \ \ {\bf L} \models \beta(x_0, x_1, \ldots, x_{n-1}) \approx 1$,

\item ${\bf L} \in {\cal V}(\mathbb{B}_{(n,3)}) \Longleftrightarrow {\bf L} \models x^\ast \wedge \gamma(y,x) \approx x^\ast \ \  \mbox{and} \ \ {\bf L} \models \beta(x_0, x_1, \ldots, x_{n-1}) \approx 1$.

\end{itemize}

\end{Theorem}

The above results give  equational bases for all join-irreducible varieties in the lattice $\Lambda({\cal BPK})$. In order to find equational bases for the remaining subvarieties, observe that $\beta(x_0,x_1, \ldots, x_{n-1}) \in [0,d] \cup \{1\}$ for each finite subdirecly irreducible $bpk$-algebra. From this, we obtain the following result.

\begin{Proposition}
Given a $pk$-algebra ${\bf L}$:

\begin{enumerate}
\item[$(1)$] ${\bf L} \in {\cal BPK}_0 \ \vee \ {\cal V}(\mathbb{B}_{(n,3)}) \Longleftrightarrow {\bf L} \models x^\ast \wedge \gamma(y,x) \approx x^\ast \ \  \mbox{and} \ \ {\bf L} \models C(y)^\prime \leq (C(y) \vee \beta(x_0, x_1, \ldots, x_{n-1}))$.

\item[$(2)$] ${\bf L} \in {\cal BPK}_1 \ \vee \ {\cal V}(\mathbb{B}_{(n,3)}) \Longleftrightarrow {\bf L} \models x^\ast \wedge \gamma(y,x) \approx x^\ast \ \  \mbox{and} \ \ {\bf L} \models C(y) \wedge C(y)^\prime \leq (C(z) \vee \beta(x_0, x_1, \ldots, x_{n-1}))  $.

\end{enumerate}
 \end{Proposition}

 \begin{Proof}
In order to prove $(1)$, note that if ${\bf L} \models C(y)^\prime \leq (C(y) \vee \beta(x_0, x_1, \ldots, x_{n-1}))$ and ${\bf L} \models \beta(x_0, x_1, \ldots, x_{n-1}) \approx 1$, then  ${\bf L} \in {\cal V}(\mathbb{B}_{(n,3)})$. On the contrary, if there exist elements in $L$ such that $\beta(x_0, x_1, \ldots, x_{n-1}) \ne 1$,  we have that  $\beta(x_0, x_1, \ldots, x_{n-1}) \leq d$ and, since $d \leq C(y)$, we obtain that ${\bf L} \in {\cal BPK}_0$. We can argue analogously to prove $(2).$
 \end{Proof}

The rest of the subvarieties can be characterized taking into account the following equalities.

\begin{itemize}

\item ${\cal V}(\mathbb{B}_{(n_1,1)}) \vee {\cal V}(\mathbb{B}_{(n_2,2)}) \vee {\cal V}(\mathbb{B}_{(n_3,3)}) ={\cal V}(\mathbb{B}_{(n_1,3)}) \cap ({\cal BPK}_1 \vee {\cal V}(\mathbb{B}_{(n_3,3)})) \cap ({\cal BPK}_0 \vee {\cal V}(\mathbb{B}_{(n_2,3)})),$

\item  ${\cal V}(\mathbb{B}_{(n_1,1)}) \vee {\cal V}(\mathbb{B}_{(n_2,2)}) = {\cal BPK}_1 \cap {\cal V}(\mathbb{B}_{(n_1,3)}) \cap ({\cal BPK}_0 \vee {\cal V}(\mathbb{B}_{(n_2,3)})),$

\item ${\cal BPK}_0 \vee {\cal V}(\mathbb{B}_{(n_2,2)}) \vee {\cal V}(\mathbb{B}_{(i_3,3)}) = ({\cal BPK}_0 \vee {\cal V}(\mathbb{B}_{(n_2,3)})) \cap ({\cal BPK}_1 \vee {\cal V}(\mathbb{B}_{(n_3,3)})),$

\item ${\cal BPK}_0\vee {\cal V}(\mathbb{B}_{(n_2,2)}) ={\cal BPK}_1 \cap ({\cal BPK}_0 \vee {\cal V}(\mathbb{B}_{(n_2,3)}))$.

\end{itemize}

\section{Free algebras on ${\cal BPK}_0$}

${\cal BPK}_0$ is the subvariety generated by simple members of ${\cal BPK}$. Next, we will determine the structure of ${\cal F}_{{\cal BPK}_0}(G)$, the free algebra over a finite set in the variety ${\cal BPK}_0$. Here and subsequently we consider $|G|=n$.

Recall that  ${\cal BPK}_0$ is the subvariety of $pk$-algebras generated by the subdirectly irreducibles  algebras of Type 1, that is, $\mathbb{B}_{(k,1)} =\mathbb{B}^-_k \oplus {\bf C}_1 \oplus \mathbb{B}^+_k$. In order to simplify the notation we denote these algebras by $\mathbb{B}_k$, that is $\mathbb{B}_k = \mathbb{B}^-_k \oplus \mathbb{B}^+_k$, where $\mathbb{B}^-_k$ and $\mathbb{B}^+_k$ have  an underlying structure of Boolean lattice with $k$ atoms.

We claim that ${\cal BPK}_0$ is a discriminator variety. In order to prove that, we observe that the subdirectly irreducible algebras in ${\cal  BPK}_0$ are regular algebras, that is, $x^\ast = y^\ast $ and $x^{\prime\ast} = y^{\prime\ast}$ implies $ x =y$. Hence, for each algebra in ${\cal BPK}_0$, it is possible to define a Heyting implication in terms of $\vee, \wedge$ and $\ast$ as follows (for more details see \cite{San}):

 $$a \to b=(a^* \vee
b^{**})^{**} \wedge [(a \vee a^*)^+ \vee a^* \vee b \vee b^*], $$ where $a^+=a^{'*'}$.

In addition, if we consider the term

$$F(x)= T(x)^* \wedge x^{**}=\left\{
                              \begin{array}{ll}
                                0 & \hbox{ if } x \ne 1 \\
                                1 & \hbox{ if } x=1
                              \end{array}
                            \right.$$
where $T(x)=C(x) \wedge C(x)^\prime = [(x \wedge x') \vee (x \wedge x')^*] \wedge [(x \wedge x') \vee  (x \wedge x')^*]'$,  it is immediate that a discriminator term for each subdirectly irreducible algebra in ${\cal BPK}_0$ is given by
$$t(x,y,z)= [F((x \to y) \wedge (y \to x)) \wedge z] \vee [(F((x \to y) \wedge (y \to x)))^* \wedge x].$$

Since ${\cal BPK}_0$ is a discriminator variety, it is known (see \cite[p.187]{BurSan})  that each non-trivial finite algebra ${\bf L}$ in the variety ${\cal BPK}_0$ is isomorphic to a  finite product of  simple algebras. In addition, since ${\cal BPK}_0$ is a locally finite variety, we have that ${\cal F}_{{\cal BPK}_0}(G)$ is finite and consequently has a factorization as

 $${\cal
F}_{{\cal BPK}_0}(G) \cong \prod_{\theta \in {\cal M}} {\mathcal
F}_{{\cal BPK}_0}(G)/\theta,$$ where ${\cal M}$
is the set of the maximal congruences on ${\mathcal
F}_{{\cal BPK}_0}(G)$.

\

Furthermore, taking into account the results of the previous section, we know that each direct factor ${\mathcal
F}_{{\cal BPK}_0}(G)/\theta$ of the free algebra in ${\cal BPK}_0$ is a simple algebra in ${\cal BPK}_0$, that is, it is isomorphic to $\mathbb{B}_k$ or ${\bf 2}$ (two-element chain).

\
Since each simple algebra ${\bf L}$ in ${\cal BPK}_0$ is an ordinal sum of Boolean algebras, the following results are easy to check.

\begin{Lemma} \label{Bk generada equivale a Boole  n generada}
For $k \geq 2$, $\mathbb{B}_{k}= \mathbb{B}_k^- \oplus \mathbb{B}_k^+$ satisfies the following:

\begin{enumerate}

\item If $G$ is a set of generators of $\mathbb{B}_k$, then there exists a set $\tilde{G} \subseteq \mathbb{B}_k^-$, $|\tilde{G}| \leq |G|$, such that  $\mathbb{B}_k$ is generated by $\tilde{G}$.

\item Given $G \subseteq \mathbb{B}_k^-$, the Boolean algebra $\mathbb{B}^-_k$ is generated by $G$ if and only if $\mathbb{B}_k$ is generated by $G$.
\end{enumerate}

\end{Lemma}

From this and from Lemma \ref{Bolle n generada equivale a Bk n generada} it follows that if $\mathbb{B}_k$ is $n$-generated then the Boolean algebra ${\mathbb B}^-_k$ is also $n$-generated. Taking into account that the largest $n$-generated  Boolean algebra is the free Boolean algebra ${\bf B}_{2^n}$, and that there exists a surjective Boolean homomorphism  $h: {\cal F}_{\cal B}(G) \to {\bf B}_k$ for each $1 \leq k\leq 2^n$ we have the following result.

\begin{Proposition}
For each $1 \leq k \leq 2^n$,  there exists a maximal congruence $\theta$ on ${\cal F}_{{\cal BPK}_0}(G)$ such that ${\cal F}_{{\cal BPK}_0}(G) / \theta \cong \mathbb{B}_k$.
\end{Proposition}

Our next goal is to decide how many  factors in the free algebra ${\cal F}_{{\cal BPK}_0}(G)$ are isomorphic to $\mathbb{B}_k$, for each $1 \leq k \leq 2^n$, that is, we want to determine $|{\cal M}_k|$, where ${\cal M}_k  =\{ \theta \in Con
(\mathcal{F}_{{\cal BPK}_0}(G)): \mathcal{F}_{{\cal BPK}_0}(G)/\theta \cong
{\mathbb B_k}\}$,  for $k \geq 1$. We also need to determine the number of congruences $\theta$ on  $\mathcal{F}_{{\cal BPK}_0}(G)$  such that $\mathcal{F}_{{\cal BPK}_0}(G)/\theta$ is isomorphic to ${\bf 2}$. Let ${\cal M}_0=\{ \theta \in Con
(\mathcal{F}_{{\cal BPK}_0}(G)): \mathcal{F}_{{\cal BPK}_0}(G)/\theta \cong
{\bf 2}\}$.

It is easily seen that

\begin{equation}
\label{cardinal de Mk} |{\cal M}_k| =\displaystyle \frac{|Sur({\cal F}_{{\cal BPK}_0}(G),
\mathbb{B}_k)|}{|Aut(\mathbb{B}_k)|},
\end{equation}
 where $Sur({\cal F}_{{\cal BPK}_0}(G),
\mathbb{B}_k)$  denotes the set of surjective homomorphisms from ${\cal F}_{{\cal BPK}_0}(G)$ to $\mathbb{B}_k$ and $Aut(\mathbb{B}_k)$ denotes the set of automorphisms on $\mathbb{B}_k$.

We now determine the number of elements of $H_k=Sur({\cal F}_{{\cal BPK}_0}(G),
\mathbb{B}_k)$ and $A_k=Aut(\mathbb{B}_k)$.

A known result about Boolean algebras states that  $|Sur({\cal F}_{\cal B}(G),
{\bf B}_k)|= \displaystyle \frac{2^n!}{(2^n-k)!}$ and  $|Aut({\bf B}_k)|=k!$ and it is immediate by Lemma \ref{Bk generada equivale a Boole  n generada} that $|Aut(\mathbb{B}_k)| = |Aut({\bf B}_k)|$. Nevertheless, $|Sur({\cal F}_{{\cal BPK}_0}(G),
\mathbb{B}_k)| \ne |Sur({\cal F}_{{\cal B}}(G),
{\bf B}_k)|$.  In order to determine $|Sur({\cal F}_{{\cal BPK}_0}(G),
\mathbb{B}_k)| $, note that surjective homomorphisms $h: {\cal F}_{{\cal BPK}_0}(G) \to
\mathbb{B}_k$ are in an obvious  one-to-one correspondence with $n$-tuples that generate $\mathbb{B}_k$  by $h \longmapsto (h(g_1), g(g_2), \ldots, h(g_n))$. In the sequel, we will say that $(x_1, \ldots, x_n)$ generates $\mathbb{B}_k$ if means $\{x_1, \ldots, x_n\}$ generates $\mathbb{B}_k$. We will denote by ${\cal T}$ the set of $n$-uples that generate $\mathbb{B}_k$. In order to calculate $|Sur({\cal F}_{{\cal BPK}_0}(G),
\mathbb{B}_k)|=|{\cal T}|$ we introduce the following definition.

\begin{Definition}
Given $\mathbb{B}_k$ and ${\bf u}=(x_1, x_2, \ldots, x_n)$ such that $x_i \in \mathbb{B}_k^-$, we define $S_{\bf u}=\{(y_1, y_2, \ldots, y_n) : y_i=x_i \mbox{ or } y_i=x_i^\prime\}$.
\end{Definition}

It is easy to check that $\{S_{\bf u}, {\bf u} \in {\cal T} \cap (\mathbb{B}_k^-)^n\}$ is a partition of ${\cal T}$. Hence, $|{\cal T}|= \displaystyle \sum_{{\bf u} \in {\cal T} \cap (\mathbb{B}^-_k)^n} |S_{\bf u}|$. Moreover, by Lemma \ref{Bk generada equivale a Boole  n generada}, $|{\cal T}|= \displaystyle \sum |S_{\bf u}|$, where ${\bf u}$ generates the Boolean algebra $\mathbb{B}_k^-$.

If ${\bf u}=(x_1, x_2, \ldots, x_n)$ generates the Boolean algebra $\mathbb{B}_k^-$ and $x_i \ne 1^-=d$ then $|S_{\bf u}|= 2^n$ (putting $x_i$ or $x_i^\prime$ in each component). This is not the case for $x_i=1^-$ because $d=d^\prime$ in $\mathbb{B}_k$.
Then, it is necessary to know how many components are equal to 1 in each $n$-tuple that generates $\mathbb{B}_k^-$ or, equivalently, given a surjective homomorphism $h:{\cal F}_{\cal B}(G) \to {\bf B}_k$, for which free generators  $g$ we have $h(g)=1$.

Each surjective homomorphism $h: {\cal F}_{\cal B}(G) \to {\bf B}_k$ has an associated  injective continuous function $f: \mathbb{X}({\bf B}_k) \to \mathbb{X}({\cal F}_{\cal B}(G))$ and each one of these is associated with a $k$-tuple in $\mathbb{X}({\cal F}_{\cal B}(G))$.

 Recall that the family of subsets of ${\bf 2}^n$ is the free Boolean algebra with $n$ free generator. Observe that ${\cal F}_{\cal B}(G)$ has $2^n$ atoms. The elements of $\mathbb{X}({\cal F}_{\cal B}(G))$  are represented by  $n$-tuples $(x_1, x_2, \ldots, x_n)$ and the clopen increasing sets
$g_i=\{(x_1, x_2, \ldots, x_n) \in {\bf 2}^n, x_i =1\}$, $1 \leq i \leq n$ are a set of  free generators.

Taking into account the dual space described above, each $k$-tuple in $\mathbb{X}({\cal F}_{\cal B}(G))$  can be represented in the following matrix

\begin{equation}
\label{matrix}
\begin{pmatrix}
  x_{11} & x_{12} & \ldots & x_{1n} \\
  x_{21} & x_{22} & \ldots & x_{2n}\\
  \vdots & \vdots & \vdots & \vdots \\
  x_{k1} & x_{k2} & \ldots & x_{kn}\\
\end{pmatrix}
\end{equation}
where each row represents an element of $\mathbb{X}({\cal F}_{\cal B}(G))$.

Let $f:\mathbb{X}({\bf B}_k) \to \mathbb{X}({\cal F}_{\cal B}(G))$ be an injective morphism and let $\mathbb{D}(f): \mathbb{D}(\mathbb{ X}({\cal F}_{\cal B}(G)))\to \mathbb{D}(\mathbb{X}({\bf B}_k))$ be its associated surjective homomorphism. $\mathbb{D}(f)(g_i)=1$, where $g_i$ is a generator, if under the duality $\mathbb{D}(f)(V_{g_i}) = f^{-1}(V_{g_i})=\mathbb{X}({\bf B}_k)$ and this is also equivalent to the  elements of column $i$  in (\ref{matrix}) all being 1.

Write $g(n,k)= |\{h \in Sur({\cal F}_{\cal B}(G),
{\bf B}_k): h(g) \ne 1, \mbox{ for each } g \in G\}|$. Observe that $g(0,k)=1$ and $g(n,k)=0$ for each $k > 2^n$.

Subtracting from $|Sur({\cal F}_{{\cal B}}(G),
{\bf B}_k)|= \frac{2^n!}{(2^n-k)!}$ the number of surjective homomorphisms that map exactly one generator to 1, then the number of surjective homomorphisms that map exactly two generators to 1, etc., we obtain the following recurrence formula.

$$g(n, k) = \frac{2^n!}{(2^n-k)!} - \displaystyle \sum_{i=1}^n \begin{pmatrix}
                                                            n \\
                                                            i \\
                                                          \end{pmatrix} g(n-i, k)$$

We use this result to find a formula to compute the number of surjective homomorphisms $h: {\cal F}_{\cal B}(G) \to {\bf B}_k$ that map exactly $i$ generators to 1, that is, whose associated matrix has exactly $i$ columns filled with ones. This number is given by:

$$f(n,k,i)= \begin{pmatrix}
              n  \\
              i \\
            \end{pmatrix} g(n-i,k).$$

\noindent The number $\begin{pmatrix}
              n  \\
              i \\
            \end{pmatrix}$ results from considering the different choices for the $i$ columns filled whith ones.

Note that given an algebra $\mathbb{B}_k$, $f(n,k,i)$ is the number of $n$-tuples $(x_1, x_2, \ldots, x_n)$ that have exactly $i$ components equal to $1^-$ and generate $\mathbb{B}_k^-$. Furthermore observe that associated to each for these $n$-tuples there are $2^{n-i}$ different $n$-tuples that generate $\mathbb{B}_k$, that is $|S_{(x_1, x_2, \ldots, x_n)}|=2^{n-i}$.  It follows immediately that

$$|Sur({\cal F}_{{\cal BPK}_0}(G)),
\mathbb{B}_k)|= \sum_{i=0}^n2^{n-i}f(n,k,i).$$

It is easy to check that $|{\cal M}_0|=2^n$ and $|{\cal M}_1|=3^n-2^n$. Finally, from this and ($\ref{cardinal de Mk}$) we have that the free  algebra over a finite poset in the variety ${\cal BPK}_0$ is given by

            $${\cal F}_{{\cal BPK}_0}(G) \cong {\bf 2}^{2^n} \times {\bf 3}^{3^n-2^n} \times \prod_{k=2}^{2^n} \mathbb{B}_k^{\frac{\sum_{i=0}^n 2^{n-i}f(n,k,i)}{k!}}$$

Using a simple computer program we can calculate the values of $f(n,k,i)$ for small $n$ and $k$, and obtain  some examples of the free algebras in the variety ${\cal BPK}_0$.

\begin{itemize}
\item ${\mathcal F}_{{\cal BPK}_0}(1) \cong  {\bf 2}^2 \times {\bf 3} \times \mathbb{B}_2^2.$

\item ${\mathcal F}_{{\cal BPK}_0}(2) \cong  {\bf 2}^4 \times {\bf 3}^5 \times \mathbb{B}_2^{20} \times \mathbb{B}_3^{16} \times \mathbb{B}_4^4$

\item ${\mathcal F}_{{\cal BPK}_0}(3) \cong  {\bf 2}^8 \times {\bf 3}^{19} \times \mathbb{B}_2^{158} \times \mathbb{B}_3^{400} \times \mathbb{B}_4^{548} \times \mathbb{B}_5^{448} \times \mathbb{B}_6^{224} \times \mathbb{B}_7^{64} \times \mathbb{B}_8^8$
\end{itemize}

\section{Finite weakly projective algebras in ${\cal BPK}_0$}

In this section we characterize the finite weakly projective algebras in the variety ${\cal BPK}_0$. Recall that an algebra ${\bf L}$ in a variety ${\cal V}$ is {\it weakly projective} if whenever $f:{\bf L}_1 \to {\bf L}_2$ is a surjective homomorphism between algebras ${\bf L}_1, {\bf L}_2 \in {\cal V}$ and $h:{\bf L} \to {\bf L}_2$ is a homomorphism, there exists a homomorphism $g:{\bf L}\to {\bf L}_1$ such that $h=f \circ g$, in other words, the following diagram commutes:

$$\begin{tikzpicture}

\node at (0,0){${\bf L}_2$};
\node at (2.1,0){${\bf L}$};
\node at (0,2){${\bf L}_1$};

\node at (-0.3,1.1){$f$};

\node at (1.2,-0.3){$h$};

\node at (1.4,1.1){$g$};

\draw[<-, >=latex]  (0.3,0) -- (1.8,0);

\draw[->, >=latex]  (0,0.5) -- (0,0.2);

\draw[->, >=latex]  (0,1.8) -- (0,0.35);

\draw[->, >=latex][dashed] (1.9, 0.2) -- (0.2, 1.8);

\end{tikzpicture}$$

\begin{Remark} {\em
Recall that in the definition of projective objects in category theory epimorphism are considered instead of surjective homomorphisms. It is known that in the category associated with any variety of algebras, mononmorphism coincide with injective homomorphisms. However, although every surjective homomorphism is an epimorphism, the converse relation does not hold in general. In ${\cal BPK}_0$, for example, if we consider the chains ${\bf 2}$ and ${\bf 3}$ and the homomorphism $f:{\bf 2} \to {\bf 3}$ given by $f(0)=0$ and $f(1)=1$, we have that $f$ is an epimorphism. In fact, if $h,g$ are homomorphisms such that $g \circ f = h \circ f$ then  necessarily $h=g$ (note that $d=d^\prime$ implies $h(d)=g(d)$). Nevertheless it is obvious that $f$ is not a surjective homomorphism.}
\end{Remark}

We recall now that given algebras ${\bf L}$ and ${\bf L}_1$ in a variety ${\cal V}$, ${\bf L}$ is a ${\it retract}$ of ${\bf L}_1$ if there exist homomorphisms $f:{\bf L} \to {\bf L}_1$ and $g: {\bf L}_1 \to {\bf L}$ such that $g \circ f = id_{\bf L}$. It is well known that ${\bf L}$ is a finite weakly projective algebra if only if ${\bf L}$ is a retract of the free algebra ${\cal F}_{\cal V}(n)$, with $n < \omega$.

\begin{Proposition} \label{proyectivas no tienen punto fijo}
If ${\bf L} \in {\cal BPK}_0$ is a finite algebra with a fixed point ($x=x^\prime$), then ${\bf L}$ is not a weakly projective algebra.
\end{Proposition}

\begin{Proof}
Let $x=x^\prime$ be a fixed point of ${\bf L}$. If we consider ${\bf L}_2 = {\bf L}$, ${\bf L}_1={\bf 2} \times {\bf L}$, the projection map $\pi_{\bf L}: {\bf L}_1 \to {\bf L}_2$, and  the identity map $id: {\bf L}_2 \to {\bf L}_2$, then there exists no homomorphism from ${\bf L}$ to ${\bf L}_1= {\bf 2} \times {\bf L}$ (${\bf L}$ has a fix point but ${\bf L}_1$ does not have one).

\end{Proof}

Let ${\bf L}$ be a finite algebra in $ {\cal BPK}_0$. Since ${\cal BPK}_0$ is a discriminator variety, ${\bf L}$ is a product of simple algebras, that is, ${\bf L} = \displaystyle \prod_{i=1}^n {\bf D}_i$, where ${\bf D}_i \cong \mathbb{B}_{k_i}$. By Proposition \ref{proyectivas no tienen punto fijo},  if ${\bf L}$ is a weakly projective algebra, there exists $k_i$ such that ${\bf D}_i \cong {\bf 2}$.  Moreover, we have the following results.

\begin{Proposition}
If ${\bf L} = {\bf 2} \times {\bf A}$ is a finite algebra in ${\cal BPK}_0$, then ${\bf L}$ is a weakly projective algebra in ${\cal BPK}_0$.
\end{Proposition}

\begin{Proof}
Let ${\bf L}={\bf 2} \times {\bf A}$ be a finite algebra in ${\cal BPK}_0$. It is enough to show that ${\bf L}$ is a retract of the free algebra ${\cal F}_{{\cal BPK}_0}(n)$ for some $n$. Since ${\bf L}$ is finite, there exists $n_0$ such that ${\bf L} \cong {\cal F}_{{\cal BPK}_0}(n_0) / \theta$. From Section 7, we know that ${\cal F}_{{\cal BPK}_0}(n_0) =  \prod_{i=0}^{2^{n_0}} \mathbb{B}_i^{m_i}$, where $\mathbb{B}_0 = {\bf 2}$ and $m_i= |{\cal M}_i|$. Since each component in ${\cal F}_{{\cal BPK}_0}(n_0)$ is simple for $i=0, \ldots, 2^{n_0}$ we have that ${\bf L}$ is isomorphic to a product of some of the components of ${\cal F}_{{\cal BPK}_0}(n_0)$, so ${\bf L} \cong \prod_{i=0}^{2^{n_0}} \mathbb{B}_i^{r_i}$ where $0 \leq r_i \leq m_i$. Moreover, we may assume that ${\bf L}={\bf 2} \times \prod_{i=1}^t {\bf D}_i^{r_i}$, where ${\bf D}_i$ is some of $\mathbb{B}_0, \mathbb{B}_1, \ldots, \mathbb{B}_{2^{n_0}}$. Thus, we can write ${\cal F}_{{\cal BPK}_0}(n_0) =  {\bf 2} \times \prod_{i=1}^t {\bf D}_i^{r_i} \times \prod_{j=1}^s {\bf E}_j$ where both  ${\bf D}_i $ and ${\bf E}_j$ are isomorphic to algebras of the type $\mathbb{B}_k$.
Let us consider now the map $g: {\bf 2} \times \prod_{i=1}^t {\bf D}_i^{r_i} \to {\bf 2} \times \prod_{i=1}^t {\bf D}_i^{r_i} \times \prod_{j=1}^s {\bf E}_j$ given by $(e,x) \longmapsto(e,x,e, \ldots, e)$ where $e \in {\bf 2}$ and $x \in \prod_{i=1}^t {\bf D}_i^{r_i}$. Clearly, $g$ is a homomorphism. In addition, if we consider the projection map $f: {\bf 2} \times \prod_{i=1}^t {\bf D}_i^{r_i} \times \prod_{j=1}^s {\bf E}_j \to {\bf 2} \times \prod_{i=1}^t {\bf D}_i^{r_i}$ we have that $f \circ g=id_{{\bf L}}$. This shows that ${\bf L}= {\bf 2} \times \prod_{i=1}^{t}{\bf D}_i^{r_i}$ is a retract of a free algebra in ${\cal BPK}_0$.

\end{Proof}

\begin{Corollary}
Let ${\bf L} \in {\cal BPK}_0$ be a finite algebra. The following conditions are equivalent.
\begin{enumerate}

\item ${\bf L}$ is a weakly projective algebra.

\item ${\bf L}$ does not have a fixed point.

\item  ${\bf L} \cong {\bf 2} \times {\bf A}$ where ${\bf A}$ is a finite algebra.

\end{enumerate}
\end{Corollary}

\section{Concluding remarks}


In this article we have studied the subvarieties of the variety ${\cal PK}$ depicted in the following diagram:

$$
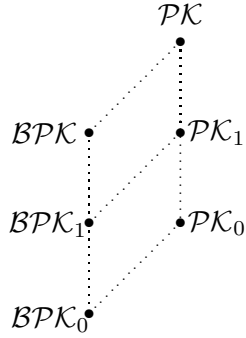
\beginpicture
\setcoordinatesystem units <1.20000cm,1.20000cm>
\setlinear
\setshadesymbol <z,z,z,z> ({\fiverm .})
\setdots <3pt>
\plot -1.00000 3.00000 -2.00000 2.00000 /
\putrule from -2.00000 2.00000 to -2.00000 1.00000
\putrule from -2.00000 1.00000 to -2.00000 0.00000
\putrule from -1.00000 3.00000 to -1.00000 2.00000
\plot -1.00000 2.00000 -2.00000 1.00000 /
\plot -1.00000 3.00000 -1.00000 1.00000 /
\plot -1.00000 1.00000 -2.00000 0.00000 /
\put {$\scriptstyle\bullet$} at -1.00000 3.00000
\put {$\scriptstyle\bullet$} at -2.00000 2.00000
\put {$\scriptstyle\bullet$} at -2.00000 1.00000
\put {$\scriptstyle\bullet$} at -2.00000 0.00000
\put {$\scriptstyle\bullet$} at -1.00000 2.00000
\put {$\scriptstyle\bullet$} at -1.00000 1.00000
\put {${\cal PK}$} at -1 3.30585
\put {${\cal BPK}$} at -2.48861 2
\put {${\cal BPK}_1$} at -2.45352 0.98201
\put {${\cal BPK}_0$} at -2.42354 -0.06747
\put {${\cal PK}_1$} at -0.60480 2
\put {${\cal PK}_0$} at -0.60 1
\endpicture$$


We have fully described the subvariety lattice of ${\cal BPK}$ as well as explicitly stated equational bases for each subvariety. From these results it follows that the pairs $({\cal V}(\mathbb{B}_{(1,2)}), {\cal BPK}_0)$ and $({\cal V}(\mathbb{B}_{(1,3)}), {\cal BPK}_1)$ are splittings pairs of the subvariety lattice of ${\cal BPK}$. Note that $\mathbb{B}_{(1,2)}$ and $\mathbb{B}_{(1,3)}$ are the 4 and 5-element chains, respectively. Thus, ${\cal BPK}_0$ and ${\cal BPK}_1$ are characterized as the largest subvarieties not containing the 4 and 5-element chains, respectively.

The problem of determining all subvarieties of ${\cal PK}$ is complex; observe that all Kleene subvarieties considered in \cite{San2020} are contained in ${\cal PK}_0$. However, two possible ways of continuing our work seem reasonable. We intend to characterize all splitting pairs in ${\cal PK}$ and to expand our knowledge of  the subvariety lattice of ${\cal PK}_0$. Since the algebras in ${\cal PK}_0$ admit a Heyting implication, we could use  the results in \cite{Cas-Mun} together with \cite{San2020} and the current article to carry this out.

\

\

\

D. Casta\~no

\noindent Departamento de Matem\'atica (Universidad Nacional del Sur) \\
Insituto de Matem\'atica (INMABB) - UNS-CONICET \\
Bah\'{i}a Blanca, Argentina

\noindent diego.castano@uns.edu.ar

\

V. Casta\~no

\noindent Departamento de Matem\'atica, Facultad de Econom\'{\i}a y Administraci\'on, Universidad Nacional del Comahue \\
Neuqu\'en, Argentina

\noindent cvaleria@gmail.com

\

J. P. D\'{\i}az Varela

\noindent Departamento de Matem\'atica (Universidad Nacional del Sur) \\
Insituto de Matem\'atica (INMABB) - UNS-CONICET \\
Bah\'{i}a Blanca, Argentina

\noindent usdiavar@criba.edu.ar

\

M. Mu\~noz Santis

\noindent Departamento de Matem\'atica, Facultad de Econom\'{\i}a y Administraci\'on, Universidad Nacional del Comahue \\
Neuqu\'en, Argentina

\noindent santis.marcela@gmail.com

\

\end{document}